\documentclass[10pt]{article}
\usepackage{geometry}                
\usepackage{graphicx}
\usepackage{subfig}
\usepackage{amssymb}
\usepackage{amsmath}
\usepackage{hyperref}
\usepackage{epstopdf}
\usepackage{caption}
\newcommand{\vecd}{\boldsymbol{\mathbf{d}}}

\newcommand{\vecf}{\boldsymbol{\mathbf{f}}}

\newcommand{\veck}{\boldsymbol{\mathbf{k}}}
\newcommand{\vecl}{\boldsymbol{\mathbf{l}}}

\newcommand{\vecw}{\boldsymbol{\mathbf{w}}}
\newcommand{\vecx}{\boldsymbol{\mathbf{x}}}

\newcommand{\matA}{\boldsymbol{\mathbf{A}}}

\newcommand{\matI}{\boldsymbol{\mathbf{I}}}

\newcommand{\matK}{\boldsymbol{\mathbf{K}}}
\newcommand{\matL}{\boldsymbol{\mathbf{L}}}

\newcommand{\matW}{\boldsymbol{\mathbf{W}}}
\newcommand{\matX}{\boldsymbol{\mathbf{X}}}

\newcommand{\matsigma}{\boldsymbol{\Sigma}}
\newcommand{\matomega}{\boldsymbol{\Omega}}
\newcommand{\matgamma}{\boldsymbol{\Gamma}}

\newcommand{\ev}{\boldsymbol{\eta}}

\newcommand{\tv}{\boldsymbol{\theta}}
\newcommand{\mv}{\boldsymbol{\mu}}

\newcommand{\bv}{\boldsymbol{\beta}}

\newcommand{\erfc}{\text{erfc}}

\title{Gaussian process regression for survival\\data with competing risks}
\author{James E. Barrett\footnote{Contact: james.j.barrett@kcl.ac.uk}\,\,\, and Anthony C. C. Coolen\\
\small{Institute of Mathematical and Molecular Biomedicine, King's College London}}
\date{Draft Version -- August 31, 2014}

\begin{document}
\maketitle

\begin{abstract}

We apply Gaussian process (GP) regression, which provides a powerful non-parametric probabilistic method of relating inputs to outputs, to survival data consisting of time-to-event and covariate measurements. In this context, the covariates are regarded as the `inputs' and the event times are the `outputs'. This allows for highly flexible inference of non-linear relationships between covariates and event times. Many existing methods, such as the ubiquitous Cox proportional hazards model, focus primarily on the hazard rate which is typically assumed to take some parametric or semi-parametric form. Our proposed model belongs to the class of accelerated failure time models where we focus on directly characterising the relationship between covariates and event times without any explicit assumptions on what form the hazard rates take. It is straightforward to include various types and combinations of censored and truncated observations. We apply our approach to both simulated and experimental data. We then apply multiple output GP regression, which can handle multiple potentially correlated outputs for each input, to competing risks survival data where multiple event types can occur. By tuning one of the model parameters we can control the extent to which the multiple outputs (the time-to-event for each risk) are dependent thus allowing the specification of correlated risks. Simulation studies suggest that in some cases assuming dependence can lead to more accurate predictions.
\end{abstract}

\section{Introduction}
\label{sec:intro}

In this work we will develop an accelerated failure time model where the event times are written as an unknown (and noise corrupted) function of the covariates. Gaussian process (GP)  regression \cite{RAS06} is used to infer the unknown function in a flexible and non-parametric manner. By specifying different \emph{kernels} in the GP prior we can probabilistically infer a wide range of qualitatively different functions. From this point of view the event times are considered `outputs' and the covariates `inputs' in a regression model. We argue that this approach is a more direct way of connecting the quantities that we have experimental access to, namely the covariates and the event times. Many existing methods of analysing survival data focus on the hazard rate. Cox's proportional hazards model \cite{COX72} is arguably the most popular such approach. These methods typically assume that the hazard rate splits into two components, one that captures the time effects and one that captures the covariate effects. Cox's model further assumes that the covariate effects are linear. It is not obvious however that this factorisation is always appropriate. Hazard rate models take a more indirect route that needs to capture both the time and covariate effects on survival outcomes whereas our approach need only capture the covariate effects, and consequently fewer assumptions are required.

The event times are transformed such that they take negative and positive values. Model parameters consist of the `noise-free' function values and these are inferred in a Bayesian manner. We compute the maximum a posteriori (MAP) solution by numerically maximising the posterior density over parameters. The hyperparameters control qualitative features of the kernel function and the overall noise level. We construct the Laplace approximation of the marginal likelihood and use that to numerically compute the MAP solution for hyperparameters.

Our model can incorporate any type of censored and truncated observations relatively easily. In addition, we obtain estimates of when the event would have occurred to individuals that were censored. We perform several simulation studies which illustrate the model's ability to infer non-monotonic relationships between the covariates and event times. We compare our model to more traditional models such as the Cox proportional hazards model, the Weibull proportional hazards model and a third model that is also based on GP regression but assumes a hazard rate similar to the Cox model but with non-linear covariate effects. We also apply our approach to experimental gene expression data.

We extend our model to the competing risks scenario by using multiple output GP regression \cite{BOYLE05}. Multiple output GP regression was originally developed for situations where multiple outputs are available corresponding to given inputs. The outputs may be statistically dependent. Again, we regard the time-to-event for different risks as the `multiple outputs' and the covariates as the `inputs'. In general, multiple output GP regression can be applied to data where each input has corresponding measurements of all (or some) of the outputs. There are two features of competing risks data that are interesting in this regard. Firstly, at most one output is available for each individual (since we only measure one event time). Secondly, once one of the outputs is observed we know that the remaining outputs must be greater than the observed output. This is because we know that remaining events would have occurred after the first reported event time. Despite these differences we will show that multiple output GP regression performs well on competing risks data.

The model can assume either independent risks or dependent risks by tuning the value of one hyperparameter. Of course, the identifiability problem \cite{TSI75} means we cannot conclude whether the risks are truly independent or not in reality. Nevertheless, within the framework of the model we will infer the value of the parameter that best explains the observed data. If the assumption of dependence has a higher probability then the model will follow this, and as we will show, exploit it to potentially make more accurate predictions. Consider, for example, two strongly dependent risks. If there is a region of the covariate space where only the first event has been observed we can still make accurate predictions of when the second type of event would occur for new individuals. This is because we know the second risk will behave similarly to the first risk. We also examine the issue of what happens in the hypothetical scenario where we `disable' or `switch off' one or more risks. The joint event time density takes a particularly convenient form in our model since the event times are conditionally independent given the underlying noise-free function values. Quantities such as the marginal survival probabilities are straightforward to compute.

Existing approaches to analysing competing risks survival data commonly assume parametric or semi-parametric cause specific hazard rates. This is useful to establish whether a certain covariate is associated with a particular risk. It may be less clear, however, how a covariate is related to overall survival probabilities in the presence of competing risks since the survival function is a function of all the cause specific hazard rates. An alternative approach is to model the cumulative incidence function, using for example a form similar to a proportional hazard model \cite{FINE99}. Shared frailty models \cite{HOU00}, random effects models \cite{VAI00}, and the concepts of pseudo-values \cite{AND10} and relative survival \cite{LAM10} are other ways to analyse competing risks data. However, all of these approaches contain some parametric or semi-parametric components. Our approach differs from these strategies since we essentially focus on modelling the joint event time density in a non-parametric fashion. Similarly to the case of a single risk this avoids imposing unnecessary structural assumptions on what form the data take.

This rest of this paper is structured as follows. In Section \ref{sec:def} we apply GP regression to survival data with a single risk and independent censoring. We will develop our GP model, outline how to infer parameters and hyperparameters, and explain how to make predictions. This model is then applied to interval censored data. In Section \ref{sec:multi} we extend our model to competing risks data. We apply our model to both experimental and simulated data and present the results in Section \ref{sec:results}. We finish with discussion in Section \ref{sec:disc}.

\section{GP regression with a single risk and independent censoring}
\label{sec:def}

We firstly define a general non-linear transformation model from which several existing models can be recovered under different assumptions. This will serve as a natural starting point for our GP regression model and offer an intuitive way to compare it to existing approaches. Survival data are $D = \{(\tau_i,\Delta_i)\}_{i=1\ldots,N}$ where $\tau_i>0$ is the time until the first event for individual $i$, the indicator variable $\Delta_i=1$ means the primary event occurred first whereas $\Delta_i=0$ indicates individual $i$ was censored, and $N$ is the total number of individuals. In addition, we acquire a vector of covariate measurements $\vecx_i\in\mathbb{R}^d$ for each individual. A general transformation model assumes
\begin{equation}
\phi(\tau_i) = f(\vecx_i) + \xi_i\quad\text{for $i=1,\ldots,N$}
\label{gpsa:eq:general_model}
\end{equation}
where $\phi$ is a monotonically increasing transformation of the event times, $f(\vecx_i)$ is some function of the covariates, and $\xi_i$ is a noise random variable with a probability density function $p_{\xi}$.

Under different assumptions of $\phi$, $f$ and $p_{\xi}$ several existing models, including our GP model, can be derived as special cases of (\ref{gpsa:eq:general_model}). For example, linear transformation models \cite{CHENG95} assume $\phi$ is unspecified and $f(\vecx) = \bv\cdot\vecx$ where $\bv$ is a vector of regression weights. Various procedures for estimating the regression parameters in such models have been proposed in \cite{FINE98} and \cite{CHEN02}. Recently \cite{LU08} considered the case where $f(\vecx)$ is an unspecified smooth function and proposed a boosting estimation method based on the marginal likelihood.

If we pick $p_{\xi}(s)=\exp(s-e^s)$ and $\phi(\tau) = \log\Lambda_0(\tau)$, where $\Lambda_0(\tau)$ is the integrated baseline hazard rate, we recover models with a hazard rate similar to Cox's model:
\begin{equation}
\pi(\tau) = \lambda_0(\tau)e^{-f(\vecx)}.
\label{eq:cox_hazard}
\end{equation}
The baseline hazard rate is $\lambda_0(\tau)$. When $f(\vecx) = -\bv\cdot\vecx$ we recover Cox's original proportional hazards model. Frailty models \cite{VAU79} can be retrieved by assuming $f(\vecx) = -\bv\cdot\vecx + w$ where $w$ is a frailty term. Generalised additive models \cite{FAH11} assume $f(\vecx) = \bv\cdot\vecx + \sum_{\mu=1}^d g_{\mu}(x_{\mu})$ where $g_{\mu}$ are non-linear functions of the covariates. See \cite{MART11} and \cite{VAN13} for recent implementations of such models. Alternatively, a GP prior can be assumed for $f(\vecx)$ as shown by \cite{SAV11} and \cite{JOE12}. Viewed in this order these models seek to accommodate increasingly complicated covariate effects through more flexible and sophisticated functions of the covariates. For completeness we note that accelerated failure time models can be recovered by assuming $\phi(\tau) = \log(\tau)$ and $f(\vecx) = \bv\cdot\vecx$. Assuming different distributions for $p_{\xi}$ results in a wide variety of accelerated failure time models (see Section 2.6 of \cite{KLEIN03}).

\subsection{The GP accelerated failure time model}
\label{sec:single}

We let $t = \phi(\tau)$ denote the transformed event times. We could choose the traditional $t = \log(\tau)$ but instead we choose
\begin{equation}
t = \phi(\tau) = \log(e^{\tau/\gamma}-1).
\label{gpsa:eq:transform}
\end{equation}
This transformation has some desirable features. Provided $\gamma < \min_i(\tau_i)$ then the transformation is effectively linear. A $t=\log(\tau)$ transformation would be non-linear and this will become particularly apparent for large $\tau$ since we may have two large values of $\tau$ that once transformed are rather similar to each each other. This may make it difficult for the model to make accurate inferences for large values of $\tau$. Since we will assume Gaussian noise the uncertainty associated with large event times will be the same as for short event times but with a non-linear transformation this is not desirable. Therefore (\ref{gpsa:eq:transform}) is preferable. The distortion due to the non-linear component of the transformation (when $\tau<\gamma$) becomes apparent only during predictions. When $t$ takes negative values they are `squashed' towards the positive half of the real line. A plot of the transformation is given in the Supporting Information.

The transformation of the output variables in GP regression has been explored by \cite{SNE04}. They examine a variety of parameterised monotonic transformations and regard any transformation parameters as hyperparameters to learn during training. Their procedure infers the most appropriate transformation such that the transformed outputs can be modelled using a Gaussian process. It might be useful to apply this method in future work.

GP regression provides a powerful non-parametric probabilistic method for relating inputs $\vecx$ to outputs $t$. It is assumed that any finite collection of the noise-free outputs $f(\vecx)$ are Gaussian distributed. For compactness we write $f_i = f(\vecx_i)$. The covariance is given by the \emph{kernel} function, $k(\vecx_i,\vecx_j) = \left<(f_i - \left<f_i\right>)(f_j - \left<f_j\right>)\right>$, which roughly tells us how `similar' $\vecx_i$ and $\vecx_j$ are. We will also require the mean function $\left<f_i\right> = m(\vecx_i)$. An excellent introduction to GP regression can be found in \cite{RAS06}. We can construct a GP regression model from (\ref{gpsa:eq:general_model}) by assuming a GP prior for the noise-free function values:
\begin{equation}
p(\vecf|\matX,\tv) = \frac{e^{-\frac{1}{2}(\vecf-\ev)\cdot\matK^{-1}(\vecf-\ev)}}{(2\pi)^{N/2}|\matK|^{1/2}}
\label{gp_prior}
\end{equation}
where $[\vecf]_i = f_i$, $[\ev]_i = \eta$ with $\eta$ constant, $[\matK]_{ij} = k(\vecx_i,\vecx_j)$ is the kernel matrix, $\tv$ is a vector of hyperparameters, and $\matX$ denotes the set of $\vecx_i$ for $i=1,\ldots,N$. In this work we have used the squared exponential kernel which is defined by $k(\vecx_i,\vecx_j) = \sigma \exp(-(\vecx_i-\vecx_j)\cdot\matL(\vecx_i-\vecx_j)/2)$ where $\sigma>0$ is a hyperparameter controlling the variance of the outputs. The matrix $\matL = \text{diag}(\vecl)$ where the components of $\vecl = (l_1^{-2},\ldots,l_d^{-2})$ are known as \emph{automatic relevance determination} (ARD) parameters and roughly tell us how important each covariate is. This is because $l_{\mu}$ defines a characteristic length scale over which the output associated with covariate $\mu$ varies. If the output varies a lot with a particular covariate then it is `important'. These hyperparameters are analogous to the weights in a linear regression model or the regression coefficients in Cox regression.

For the noise variable in (\ref{gpsa:eq:general_model}) we pick $p(\xi) = \mathcal{G}(0,\beta^2)$, where $\mathcal{G}(0,\beta^2)$ denotes a Gaussian distribution with mean 0 and variance $\beta^2$. It follows that the event time density for individual $i$ is
\begin{equation}
p(t_i|f(\vecx_i)) = \mathcal{G}(f(\vecx_i),\beta^2).
\end{equation}
This has a convenient form since the conditional event time density has a simple form with all of the non linear covariate effects captured by $p(\vecf|\matX,\tv)$. From this we can derive the survival function and hazard rate:
\begin{equation}
S(\tau) = \int_{\tau}^{\infty}\text{d}s\,p(s|f_i)\quad\text{and}\quad \pi_i(\tau) = \frac{p(\tau|f_i)}{\int_{\tau}^{\infty}\text{d}s\,p(s|f_i)}.
\end{equation}

For the present section we will consider only right censoring. Interval censoring will be considered in Section \ref{sec:interval}. We infer the function values $\vecf\in\mathbb{R}^N$ using Bayes' theorem:
\begin{equation}
p(\vecf|\matX,D,\tv)= \frac{p(D|\vecf,\tv)p(\vecf|\matX,\tv)}{\int\text{d}\vecf'\,p(D|\vecf',\tv)p(\vecf'|\matX,\tv)}
\label{gpsa:eq:bayes}
\end{equation}
with $p(D|\vecf, \tv) = \prod_{i=1}^N P(t_i,\Delta_i|f_i)$ where $P(t_i,\Delta_i|f_i)$ is the likelihood contribution made by individual $i$ and depends on what type of censoring or truncation has occurred \cite[Section 3.5]{KLEIN03}. In this case non-censored individuals contribute with the event time density evaluated at the reported event time, $P(t_i,\Delta_i=1|f_i) = p(t_i|f_i)$, and a censored individual contributes the probability that the event occurred after the reported event time, $P(t_i,\Delta_i=0|f_i) = S(t_i|f_i)$. We determine the maximum a posteriori (MAP) solution by numerically minimising the negative log likelihood $\mathcal{L}(\vecf) = - N^{-1}\log p(\vecf|\matX,D,\tv)$:
\begin{equation}
\mathcal{L}(\vecf) =-\frac{1}{N}\sum_{i:\Delta_i=1}\log p(t_i|f_i) -\frac{1}{N}\sum_{i:\Delta_i=0}\log S(t_i|f_i) +\frac{1}{2N}(\vecf-\ev)\cdot\matK^{-1}(\vecf-\ev) + \frac{1}{2N}\log|\matK|.
\label{gpsa:eq:LL}
\end{equation}
Numerical optimisation is performed using a gradient based optimiser in Matlab. Partial derivatives are given in the Supporting Information. Hyperparameters are determined by optimising the Laplace approximation of the marginal likelihood $\int\text{d}\vecf'\,p(D|\vecf',\tv)p(\vecf'|\matX,\tv)$. We do this by firstly expanding $\mathcal{L}(\vecf)$ to second order around the minimum $\hat{\vecf}$ using a Taylor expansion $\mathcal{L}(\vecf) \approx \mathcal{L}(\hat{\vecf}) + \tfrac{1}{2}(\vecf-\hat{\vecf})\cdot\matA(\vecf-\hat{\vecf})$. The matrix $\matA$ contains second order partial derivatives and is defined by $\matA_{ij} = \partial^2/{\partial w_i\partial w_j}\mathcal{L}(\vecf)|_{\hat{\vecf}}$. We now write the marginal likelihood as
\begin{align}
p(D|\tv) &\approx\int\text{d}\vecw\, e^{-N\mathcal{L}(\hat{\vecf}) - \frac{N}{2}(\vecf-\hat{\vecf})\cdot\matA(\vecf-\hat{\vecf})}\nonumber\\
&= p(D|\hat{\vecf},\tv)p(\hat{\vecf}|\tv)(2\pi)^{N/2}|(N\matA)^{-1}(\tv)|^{1/2}.
\label{intro:eq:laplace}
\end{align}
We take the negative log of this
\begin{equation}
\mathcal{L}_{hyp}(\tv) = \mathcal{L}(\hat{\vecf})-\frac{1}{2}\log2\pi+\frac{1}{2N}\log|\matW+\matK^{-1}| 
\label{gpsa:eq:LLhyp}
\end{equation}
where the diagonal matrix is defined by $\matW_{ii} = -{\partial^2}/{\partial f_i^2}\log p(D|\vecf,\tv)$ (given in the Supporting Information) and numerically minimise with respect to $\tv$. Note that each evaluation of the negative hyperparameter log likelihood requires determining $\hat{\vecf}$.

\subsubsection{Predictions, hazard rates and survival curves}
\label{gpsa:sec:pred}

Having trained a GP regression model (by inferring the function values $\vecf$ and hyperparameters $\tv$) we may wish to predict the event time $\tau^*$ for a new individual with covariates $\vecx^*$. The predictive distribution for a test output $\vecf^*$ corresponding to a test input $\vecx^*$ is Gaussian with mean and variance
\begin{align}
\hat{\mu} & =\veck^*\cdot\matK^{-1}\hat{\vecf}\label{gpsa:eq:pred_mean}\\
\hat{\kappa} & = k(\vecx^*,\vecx^*) -\veck^*\cdot(\matK+\matW^{-1})^{-1}\veck^*
\label{gpsa:eq:pred_var}
\end{align}
where $[\veck^*]_i = k(\vecx^*,\vecx_i)$. These expression are similar to the usual GP predictive mean and variance except in this case we include additional variance due to the uncertainty in $\hat{\vecf}$ (see Section 3.4.2 of \cite{RAS06}). The corresponding density for (the noisy prediction) $t^*$ is $\mathcal{G}(\hat{\mu},\hat{\kappa}+\beta^2)$. Finally, we need to transform back to the original time variable $\tau^*$:
\begin{equation}
p(\tau^*|\vecx^*,\matX,D) = \frac{e^{-\frac{1}{2(\hat{\kappa}+\beta^2)}(\log(e^{\tau^*/\gamma}-1) - \hat{\mu})^2}}{(2\pi(\hat{\kappa}+\beta^2))^{1/2}}\frac{e^{\tau^*/\gamma}}{\gamma(e^{\tau^*/\gamma}-1)}.
\label{eq:predictive}
\end{equation}
Once the predictive event time density has been obtained we can compute the primary hazard rate and survival function if desired. It may also be desirable to make a specific prediction of when the event will occur. This can be done by numerically computing the mean of the event time density:
\begin{equation}
\left<\tau^*\right> = \int_0^{\infty}\text{d}s\,sp(s|\vecx^*,\matX,D).
\end{equation}
The variance $\left<(\tau^*)^2\right> - \left<\tau^*\right>^2$ can also be computed as gives us a measure of uncertainty regarding our prediction.

\subsection{GP regression with a single risk and independent interval censoring}
\label{sec:interval}

Existing methods for interval censored survival data are typically based on parametric or semi-parametric models. The advantage of parametric models is that expressions for the survival function can be obtained in closed form and hence the exact likelihood can be constructed for right, left or interval censored observations. See \cite{LIN98} for a discussion and comparison of several parametric models. Weibull accelerated failure time models are considered in \cite{OD92, RAB95, KOM09}. A family of parametric models that can handle time dependent covariates is presented in \cite{SPA06}. Most semi-parametric models are adaptations of Cox's model. The partial likelihood argument cannot be used so usually the full likelihood is numerically optimised with respect to parameters. Some representative examples can be found in \cite{SIN99, SAT96, GOE00}. Another strategy is to impute the event times \cite{LAW92} by taking the midpoint or the end of the interval for instance \cite{PAN02}, and then applying standard methods to the imputed event times.

Our GP model can readily be extended to accommodate interval censored observations. Now $\Delta=1$ corresponds to an interval censored observation and $\Delta=0$ a right censored one. We observe upper and lower times that define an interval\footnote{Note that we are still working with the transformed event times $t=\phi(\tau)$ defined by (\ref{gpsa:eq:transform}).} $(t_i^l,t_i^u)$ and we have $P(t_i^l,t_i^u,\Delta_i=1|f_i) = S(t_i^l|f_i) - S(t_i^u|f_i)$. Taking the negative log of the posterior (\ref{gpsa:eq:bayes}) and ignoring terms independent of $\vecf$ we get
\begin{align}
\mathcal{L}(\vecf) &= -\frac{1}{N}\sum_{i:\Delta_i=1}\log[S(t_i^l|f_i) - S(t_i^u|f_i)] -\frac{1}{N}\sum_{i:\Delta_i=0}\log S(t_i|f_i) -\frac{1}{N}\log p(\vecf|\matX,\tv).
\label{gpsa:eq:interval}
\end{align}
As above, we find $\hat{\vecf}$ by numerically minimising the negative log likelihood. Hyperparameters are determined using the Laplace approximation of the marginal likelihood (\ref{gpsa:eq:LLhyp}) but with a different matrix $\matW$. Inference and predictions proceed as above.

\subsection{Weibull proportional hazards model (WPHM)}
\label{sec:wphm}

For the purposes of comparison we will use a Weibull proportional hazards model. This model assumes a more traditional hazard rate (\ref{eq:cox_hazard}) with a baseline hazard rate $\lambda_0(\tau) = (\nu/\rho)(\tau/\rho)^{\nu-1}$ where $\rho>0$ is a scale parameter and $\nu>0$ is a shape parameter. The cumulative base hazard rate is $\Lambda_0(\tau) = (\tau/\rho)^{\nu}$. We infer the optimal parameter values by minimising the negative log data likelihood
\begin{equation}
\mathcal{L}(\bv,\rho,\nu)  = -\frac{1}{N}\sum_{i:\Delta_i=1}\left[\log\lambda_0(\tau_i)+\bv\cdot\vecx_i\right] +\frac{1}{N}\sum_{i=1}^N\Lambda_0(\tau_i)e^{\bv\cdot\vecx_i}.
\end{equation}
Predictions for new individuals with covariates $\vecx^*$ can be made by computing the mean (and variance) of the event time density (using optimal parameters $\hat{\bv},\hat{\rho},\hat{\nu}$)
\begin{equation}
\left<\tau\right> = \int_0^{\infty}\text{d}s\,s\lambda_0(s)e^{\hat{\bv}\cdot\vecx^*}\exp(-\Lambda_0(s)e^{\hat{\bv}\cdot\vecx^*}).
\label{eq:etd}
\end{equation}

\section{Multiple output GP regression with competing risks and independent censoring}
\label{sec:multi}

We extend the transformation model from the case of a single risk (\ref{gpsa:eq:general_model}) to the case of two competing risks (this can easily be generalised to more than two):
\begin{equation}
\phi(\tau_i^1) = f_1(\vecx_i) + \xi_i^1\quad\text{and}\quad \phi(\tau_i^2) = f_2(\vecx_i) + \xi_i^2 \quad\text{for $i=1,\ldots,N$.}
\label{multi:eq:general_model}
\end{equation}
Each event time is related to the same covariates via two different functions corrupted with two different noise random variables. In the case of competing risks the event times may be correlated so we will use multiple output GP regression to capture dependency between outputs. Multiple output GP regression was first introduced to the machine learning community by \cite{BOYLE05} who built on work developed in \cite{HIG02} which illustrated that a Gaussian process can be obtained from a convolution of a Gaussian white noise process. We will follow their approach closely in this section. The noise-free outputs are written as $f_1(\vecx) = u_1(\vecx) + s_1(\vecx)$ and $f_2(\vecx) = u_2(\vecx) + s_2(\vecx)$ where $u_1$ and $u_2$ are GPs unique to each output and $s_1$ and $s_2$ are `shared' GPs obtained by convolving the same Gaussian white noise process. Dependency between outputs can be captured via the shared  components.

The covariance between the noiseless outputs is (terms such as $\left<u_r(\vecx_i),s_{q}(\vecx_j)\right>$ vanish)
\begin{equation}
\left<f_r(\vecx_i),f_q(\vecx_j)\right> = \left<u_r(\vecx_i),u_q(\vecx_j)\right> + \left<s_r(\vecx_i),s_q(\vecx_j)\right>.
\label{eq:multi:cov}
\end{equation}
Following the example of \cite{BOYLE05} we have
\begin{equation}
\begin{array}{ll}
\left<u_r(\vecx_i),u_r(\vecx_j)\right> = \frac{\pi^{d/2}\sigma^2}{\sqrt{|\matsigma_r|}}e^{-\frac{1}{4}\vecd\cdot\matsigma_r\vecd}
& \left<s_1(\vecx_i),s_2(\vecx_j)\right> = \frac{(2\pi)^{d/2}\omega^2}{\sqrt{|\matomega_1+\matomega_2|}}e^{-\frac{1}{2}(\vecd-\mv)\cdot\matgamma(\vecd-\mv)}\\
\left<s_2(\vecx_i),s_1(\vecx_j)\right> = \frac{(2\pi)^{d/2}\omega^2}{\sqrt{|\matomega_1+\matomega_2|}}e^{-\frac{1}{2}(\vecd+\mv)\cdot\matgamma(\vecd+\mv)}
& \left<s_r(\vecx_i),s_r(\vecx_j)\right> = \frac{\pi^{d/2}\omega^2}{\sqrt{|\matomega_r|}}e^{-\frac{1}{4}\vecd\cdot\matomega_r\vecd}\\
\end{array} 
\end{equation}
where the covariance matrices are diagonal with $\matsigma_r = \Sigma_r\matI_{d\times d}$ and $\matomega_r = \Omega_r\matI_{d\times d}$ and $\matgamma = \matomega_1(\matomega_1+\matomega_2)^{-1}\matomega_2$. We assume that the characteristic length scales in covariate space are the same $\Sigma_1=\Sigma_2=\Omega_1=\Omega_2 = l^{-2}$. The vector $\mv$ allows one output to be a translation of the other. We assume $[\mu]_{\nu} = \mu$ for $\nu=1,\ldots,d$ with $\mu$ constant. Finally, we shall assume that the noise levels are the same for both events $\beta_1=\beta_2=\beta$. In the simplest case we have a six-dimensional vector of hyperparameters $\tv = (\eta,\mu,\beta,\sigma,\omega,l)$ where $\eta\in\mathbb{R}$ (from the GP mean), $\mu\in\mathbb{R},\beta\geq0,\sigma\geq0,\omega\geq0$ and $l\geq0$. These simplifications are by no means necessary and may not be appropriate for certain datasets. They do however make inference of hyperparameters considerably easier since the search space will in general contain local minima so lowering the dimension of the search space will have significant computational advantages.

Inserting these into (\ref{eq:multi:cov}) allows us to construct a covariance matrix which we can use to define a GP prior over $\vecf = [\vecf_1,\vecf_2]\in\mathbb{R}^{2N}$:
\begin{equation}
p(\vecf|\matX) = \mathcal{G}\left(\ev,
\left[
\begin{array}{cc}
\matK_{11} & \matK_{12}\\
\matK_{21} & \matK_{22}\\
\end{array}
\right]
\right)
\label{eq:multi:gpprior}
\end{equation}
where $[\matK_{rq}]_{ij} = \left<f_r(\vecx_i),f_q(\vecx_j)\right>$ and $[\ev]_{\nu} = \eta$, with $\eta\in\mathbb{R}$, is the GP mean. The block matrices have an intuitive interpretation. $\matK_{11}$ and $\matK_{22}$ control the covariance structure of the independent parts of each output whereas the off-diagonal blocks control the covariance between outputs.

Returning to (\ref{multi:eq:general_model}) we will now assume the GP prior (\ref{eq:multi:gpprior}) for the function values $\vecf = [\vecf_1,\vecf_2]$. Again we let $t_1 = \phi(\tau_1)$ and $t_2 = \phi(\tau_2)$. The indicator variable can take values $\Delta_i=0,1,2$ to indicate censoring, event type 1 or event type 2 respectively. With a Gaussian distribution for the noise we obtain $t_r \sim \mathcal{G}(f_r,\beta_r^2)$ for $r=1,2$. Assuming that right censoring is independent the joint event time density is conditionally independent given the noise-free function values
\begin{equation}
p(t_i^1,t_i^2|f_i^1,f_i^2) = p(t_i^1|f_i^1)p(t_i^2|f_i^2).
\label{multi:eq:condindep}
\end{equation}
The conditional independence leaves us with a rather convenient event time density. All of the complicated business of correlations between risks and similarities between individuals is captured by the GP prior leaving a simple product of univariate Gaussian densities. We will discus thus more in Section \ref{multi:sec:disabling}.

\subsection{Inference of noise-free function values and hyperparameters}
\label{multi:sec:inference}
We can write the data likelihood as a product of Gaussian density terms and cumulative Gaussian terms
\begin{equation}
p(D|\vecf_1,\vecf_2) =\prod_{r=1}^2\left\{\prod_{i=1}^N[p(t_i|f_i^r)]^{\delta_{\Delta_i,r}}[S(t_i|f_i^r)]^{1-\delta_{\Delta_i,r}}\right\}.
\label{multi:eq:datalikelihood}
\end{equation}
As before, we use Bayes' theorem (\ref{gpsa:eq:bayes}) to calculate the posterior over the function values and then use this to obtain the negative log likelihood
\begin{align}
\mathcal{L}(\vecf) &= -\frac{1}{N}\sum_{i:\Delta_i\neq1}\log S(t_i|f_i^1) -\frac{1}{N}\sum_{i:\Delta_i\neq2}\log S(t_i|f_i^2) -\frac{1}{N}\sum_{i:\Delta_i=1}\log p(t_i|f_i^1)\nonumber\\
& \quad\quad-\frac{1}{N}\sum_{i:\Delta_i=2}\log p(t_i|f_i^2) + \frac{1}{2N}(\vecf-\ev)\cdot\matK^{-1}(\vecf-\ev) + \log2\pi +\frac{1}{2N}\log|\matK|.
\end{align}
Hyperparameters are obtained by minimising the negative log of the Laplace approximation of the marginal likelihood. Details are given in the Supporting Information.
\subsection{Making predictions}

The predictive distribution for the output $f_*^r$ corresponding to a new input $\vecx^*$ is Gaussian with mean and variance
\begin{align}
\hat{\mu}_r &= \veck^*_r\cdot\matK^{-1}\vecf\label{multi:eq:fmean}\\
\hat{\kappa}_r &= k(\vecx^*,\vecx^*) - \veck^*_r\cdot(\matK+\matW^{-1})^{-1}\veck^*_r\label{multi:eq:fvar}
\end{align}
where $\vecf = [\vecf_1,\vecf_2]$ and $\veck^*_r =[\veck_{r1}^*,\veck_{r2}^*]$ with $[\veck_{rq}^*]_i = \left<f^r(\vecx^*),f^q(\vecx_i)\right>$ given by (\ref{eq:multi:cov}). The matrix $\matK$ is the covariance matrix in (\ref{eq:multi:gpprior}) formed out of four block matrices. Finally, $k(\vecx_*,\vecx_*) = \left<f^r(\vecx_*),f^r(\vecx_*)\right> = \pi^{d/2}\sigma^2/\sqrt{|\matsigma_r|} +\pi^{d/2}\omega^2/\sqrt{|\matomega_r|}$. The predictive density over the original event time variable is given by (\ref{eq:predictive}). From this the mean and variance can be numerically computed, similarly to Section \ref{gpsa:sec:pred}. Once the predictive event time density has been obtained one can readily derive hazard rates or survival curves for each risk if desired.

\subsection{`Disabling' a risk}
\label{multi:sec:disabling}

A perennial question in survival analysis is how to estimate the survival probabilities in the absence of one or more risks. It is not primarily a statistical questions since `disabling' or eliminating one or more risks will in general alter the remaining risks because the risks will in general share biological pathways or rely on the same biological systems. The quantities we infer from data correspond to a world where all of the risks are operating so we must assume that that these quantities are relevant to the hypothetical world where one or more risk have been somehow `disabled'. Suppose now we have a total of $R$ risks. By `disabling' all risks except risk $r$ we mean replacing
\begin{equation}
p_i(\tau_0,\ldots,\tau_R) = \tilde{p}_i^r(\tau_r)\lim_{\zeta\to\infty}\prod_{q\neq r}\delta(\tau_q-\zeta)
\label{multi:eq:marginalp}
\end{equation}
where $\tilde{p}_i^r(\tau) = \int_0^{\infty}(\prod_{q\neq r}\text{d}s_q)p(s_0,\ldots,s_R)$ is the marginal density of event time $r$. We use tildes to denote quantities after the risks have been disabled. The survival function is
\begin{equation}
\tilde{S}_i^r(\tau) = \int_{\tau_r}^{\infty}\text{d}s\,\tilde{p}_i^r(s).
\label{multi:eq:disable_survival}
\end{equation}
Since there is only one risk the cause specific hazard rate is:
\begin{equation}
\tilde{\pi}_i^r(\tau) = \tilde{p}_i^r(\tau)/\tilde{S}_i^r(\tau).
\label{multi:eq:disabled_dependents}
\end{equation}
From this it follows that 
\begin{equation}
\tilde{S}_i^r(\tau) = e^{-\int_0^{\tau}\text{d}s\,\tilde{\pi}_i^r(s)}.
\label{multi:eq:disabled_hazards}
\end{equation}
Note that $\tilde{\pi}_i^q(\tau) = 0$ and $\tilde{S}_i^q(\tau) = 1$ for all $q\neq r$. The cumulative incidence function for risk $r$ becomes
\begin{equation}
\tilde{C}_i^r(\tau) = \int_0^{\tau}  \text{d}s\,\tilde{\pi}_i^r(s)e^{-\int_0^s\text{d}s'\tilde{\pi}_i^r(s')} = 1 - \tilde{S}_i^r(\tau).
\end{equation}
The interpretation of the marginal survival probabilities depends on whether the risks are independent or not. We examine both cases separately.
\subsubsection*{Dependent risks:}
In this case $\tilde{\pi}_i^r(\tau) \neq \pi_i^r(\tau)$. This can be seen by comparing the expression for $\pi_i^r(\tau)$
\begin{equation}
\pi_i^r(\tau) = \frac{\left(\prod_{q\neq r}\int_{\tau}^{\infty}\text{d}s_q\right)p_i(s_0,\ldots,s_{r-1},\tau,s_{r+1},\ldots,s_R)}{S_i(\tau)}
\label{intro:eq:hr}
\end{equation}
to the expression for $\tilde{\pi}_i^r(\tau)$ which is given by (\ref{multi:eq:disabled_dependents}). This is to be expected since switching off the other risks will change the probability to survive until a certain time and hence the hazard rate due to risk $r$ will also change. In this case the quantity $S_i^r(\tau) = \exp(-\int_0^{\tau}\text{d}s\,\pi_i^r(s))$ cannot be interpreted as a marginal survival probability in the hypothetical world where all other risks are switched off. Consequently, $C_i^r(\tau) = 1 - S_i^r(\tau)$ does not have a valid interpretation as a cumulative probability distribution either.

\subsubsection*{Independent risks:}
In the case of independent risks the survival function can be written as $S_i(\tau) = S_i^1(\tau)\cdots S_i^R(\tau)$ where the marginal survival functions are defined as  $S_i^r(\tau) = \int_{\tau}^{\infty}\text{d}s\,p_i^r(s)$ for $r=1,\dots,R$. Since the risks are independent it immediately follows that $\tilde{p}_i^r(\tau) = p_i^r(\tau)$. From (\ref{multi:eq:disable_survival}) and (\ref{multi:eq:disabled_dependents}) it follows that $\tilde{S}_i^r(\tau) = S_i^r(\tau)$ and $\tilde{\pi}_i^r(\tau) = \pi_i^r(\tau)$. In this case the quantity $S_i^r(\tau) = \exp(-\int_0^{\tau}\text{d}s\,\pi_i^r(s))$ is equal to (\ref{multi:eq:disabled_hazards}) and hence it can be interpreted as a marginal survival probability in the hypothetical world where all other risks are switched off. 

\subsubsection*{The GP model:}

In our case the conditional independence of the event times given the noise-free function means that we can always interpret $S_i^r(\tau) = \int_{\tau}^{\infty}\text{d}s\,p_i^r(s|f_i^r)$ as a marginal survival probability. This true regardless of whether the underlying functions are independent or otherwise (which in the language of our model means this is true for any value of $\omega$).

\section{Results}
\label{sec:results}

Here we present result from simulation studies and experimental data. We begin by explaining how the simulated data are generated. We then apply our GP regression method to simulated data with a single risk and independent right censoring. We also test the performance of the model on interval censored data. Then we apply the model to gene expression data from a study of lymphoma patients \cite{ROSE02}. Finally, we generate simulated competing risks data.

\subsection{Generation of simulated data with a specified event time density}

Generation of simulated data is straightforward in the case of the GP regression model. $N$ covariate vectors are randomly generated from a uniform distribution on a finite region of the covariate space where $N$ is the number of samples we wish to generate. The corresponding kernel matrix $\matK$ is constructed, and event times are sampled from the GP prior (\ref{gp_prior}) which in practice means drawing a random vector $\vecf$ from an $N$-dimensional multivariate Gaussian density, and then adding Gaussian noise to the components of $\vecf$. Finally independent right censoring is simulated by randomly selecting a subset of the individuals and generating a random number from a uniform distribution defined on the interval $[0,\tau_i)$ which is then recorded as the time of censoring. Competing risks data are generated in the the same way but with the multiple output GP prior (\ref{eq:multi:gpprior}).

\subsection{Non-monotonic simulated survival data with a single risk}

\begin{figure}[h!]
\centering
\includegraphics[scale=0.61]{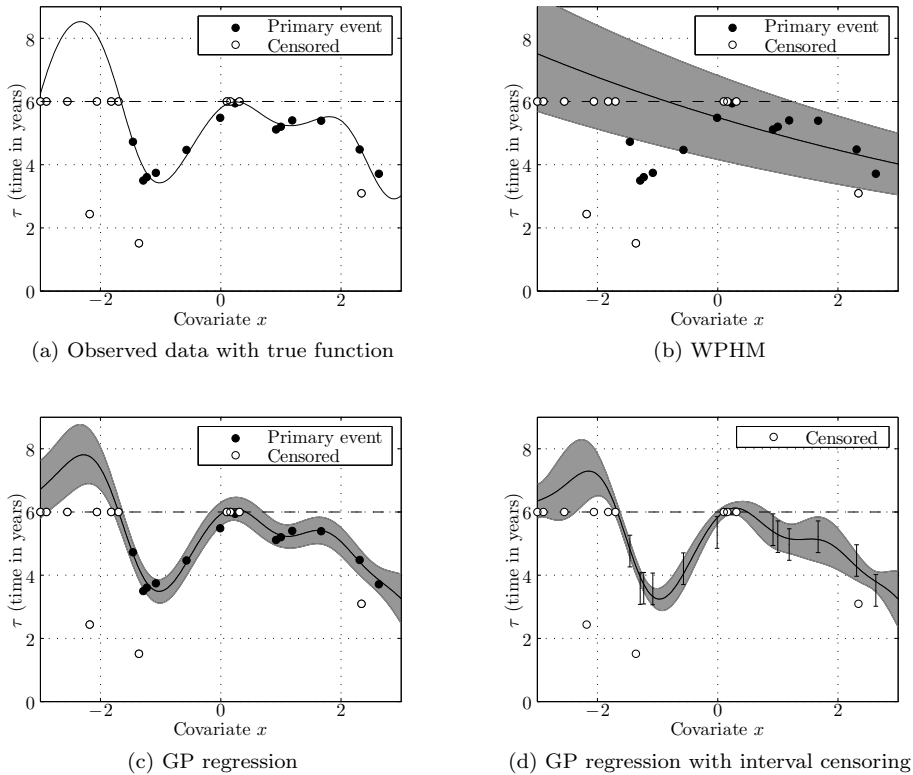}
\caption[Example of GP regression on non-monotonic simulated data]{Results from a simulated dataset with $d=1$ generated with a squared exponential kernel with hyperparameters set to $(\eta,\beta,\sigma,l) = (5,0.2,3,0.7)$. There are $N=25$ individuals, 13 of which are censored. The end of trail at 6 years is represented by the dashed line. Figure (a) shows the observed data with the `true' function. Figure (b) is a plot of the predicted event time using the Weibull proportional hazards model. The grey region represents plus and minus one standard deviation from the mean prediction. We found $(\beta,\rho,\nu) = (0.49,6.0, 4.7)$. The WPHM fails to infer the correct function since it assumes a monotonic relationship between covariates and event times. In (c) the mean prediction using our model is shown (i.e. we have plotted (\ref{eq:predictive}) as a function of $x$). Optimal hyperparameters were found to be $(\eta,\beta,\sigma,l) =  (5.82,0.32,2.59,0.64)$. Note the increased uncertainty at $x\in(-3,-2)$ where only censored observations were made. In (d) the non-censored observations were converted to interval censored observations by generating a random one year interval which are represented by the error bars. The inferred hyperparameters are $(\eta,\beta,\sigma,l) =  (5.67,0.14,3.34,0.57)$. The GP model is clearly capable of recovering non-monotonic relationships.}
\label{gpsa:fig:example1} 
\end{figure}

Shown in Figure \ref{gpsa:fig:example1} are results from a simulated dataset that consists of $N=25$ individuals with a single covariate $x$. There are 13 censored individuals and 12 who have experienced the primary risk. An end of trail cutoff at 6 years has been imposed and several individuals have been censored due to this (see Figure \ref{gpsa:fig:example1} (a)).

In Figure \ref{gpsa:fig:example1} (b) we have plotted the predicted mean event time using the Weibull proportional hazards model. The WPHM is poorly suited to these data as it assumes a monotonically increasing or decreasing relationship between event times and covariates. The results from our model are shown in Figure \ref{gpsa:fig:example1} (c). The model infers the underlying function and retrieves the hyperparameters reasonably well. The inferred function gives an estimate of when event times will occur. Note that the model has extrapolated the underlying function beyond the end of trial cutoff. This can be seen in the region $x\in(-3,-2)$, and the uncertainty is also greatest in this region. In Figure \ref{gpsa:fig:example1} (d) we convert these data into interval censored data by generating a random one year interval for all of the non-censored individuals. These intervals are represented by the `error bars' in the plot. The GP regression model is capable of recovering the underlying function.

We also implemented the GP hazard rate model from \cite{JOE12}. Additions results are available in the Supporting Information and show that the GP hazard rate model is also capable of inferring non-linear relationships and offers comparable performance to our GP model. A disadvantage with the hazard rate model is the difficulty in interpreting the hyperparameters. This is because the function inferred in that case describes the relationship between the covariates and the hazard rates. Hyperparameters such as the noise level and overall variance have a less intuitive interpretation.

\subsection{Experimental gene expression data}

\begin{figure}
\centering
\begin{tabular}{c c}
\subfloat[GP regression training data]{\includegraphics[scale=0.7]{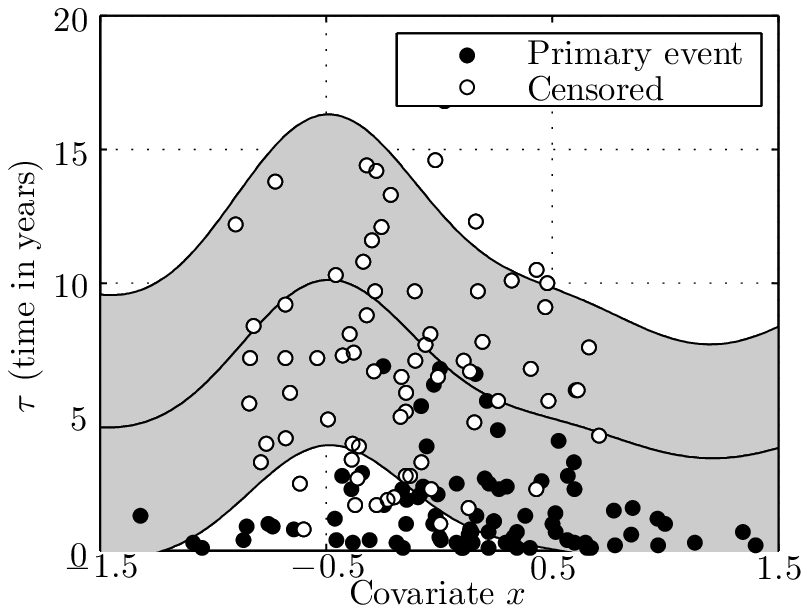}} & \subfloat[GP regression validation data]{\includegraphics[scale = 0.7]{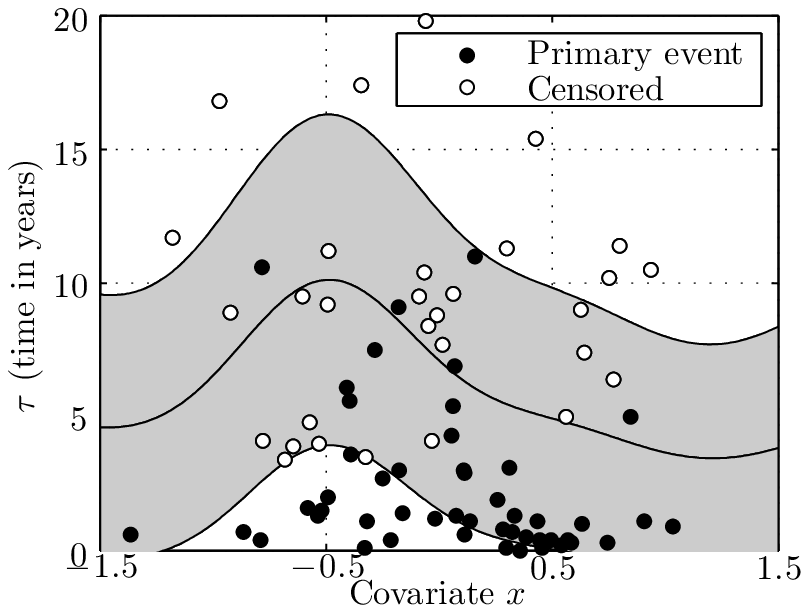}}\\
\end{tabular}
\caption[Application of GP regression to experimental gene expression data]{Univariate analysis of gene number 33014 from the Rosenwald lymphoma dataset \cite{ROSE02}. In (a) is the function inferred on the training set (of 160 patients) using our GP regression method which clearly shows a non-linear relationship between the expression levels and event times. In (b) is the same function superimposed on the validation set.}
\label{gpsa:fig:rosenwald}
\end{figure}

We applied our method to the gene expression data from the Rosenwald 2002 study of lymphoma patients \cite{ROSE02}. These data consist of $N=240$ patients each with $d=7399$ gene expression measurements. In the original analysis the patients had been split into a training group of 160 and a validation group of 80 individuals. These data were studied by \cite{LU08} to test a transformation model with non-linear covariate effects and it was reported that some of the gene expression levels had a non-linear relationship with the time-to-event. We examined one of these genes, with UNIQID $=33014$, with our GP method and also found a non-linear function $f(\vecx)$. This function was inferred using the 160 training individuals and can be seen in Figure \ref{gpsa:fig:rosenwald} (a). If we compare this to the top right panel of Figure 2 in \cite{LU08} we can see that both functions are very similar (once we ignore the fact that, by definition, they differ in sign).

To further quantify the difference between our GP method and the WPHM we computed the mean square error (MSE) between the predicted time-to-event and the reported time-to-event in both the training and validation sets for both models. The results are displayed in Table \ref{multi:tab:results2}. It is clear that the GP method offers vastly superior performance compared the WPHM. We can also see that the WPHM validation error is considerably larger than the training error. This is a hallmark of overfitting where the model fails to generalise well to unseen data. GP regression on the other hand does not suffer from this problem on this dataset.

\begin{table}[h]
  \begin{center}
    \begin{tabular}{|c|c|c|}
    \hline
     & GP regression & WPHM \\
	\hline
	Training MSE (years$^2$)& 22.86 & 774.96\\
    Validation MSE (years$^2$) & 22.38 & 2514.7\\
    \hline
    \end{tabular}
  \end{center}
\caption[GP regression applied to experimental data]{Comparison of mean square error (MSE) between the predicted and reported event times in the validation set using gene number 33014 from the Rosenwald lymphoma dataset \cite{ROSE02}. Our GP regression method offers superior performance to the WPHM because it has detected a non-linear relationship between the event times and that gene expression level. The WPHM also overfits the validation data since the MSE is considerably larger than the training MSE. The GP model does not suffer from this problem in this case.}
\label{multi:tab:results2}
\end{table}

\subsection{Comparison of GP models with dependent and independent competing risks}
\label{multi:sec:comp}

\begin{figure}[htbp!]
\centering
\begin{tabular}{c c}
\subfloat[Dependent risk 1]{\includegraphics[scale=0.7]{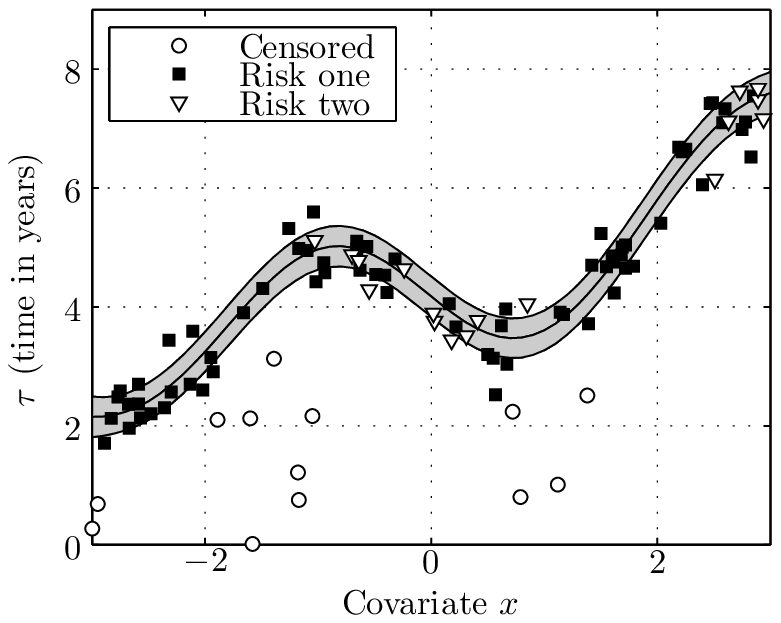}} & \subfloat[Dependent risk 2]{\includegraphics[scale = 0.7]{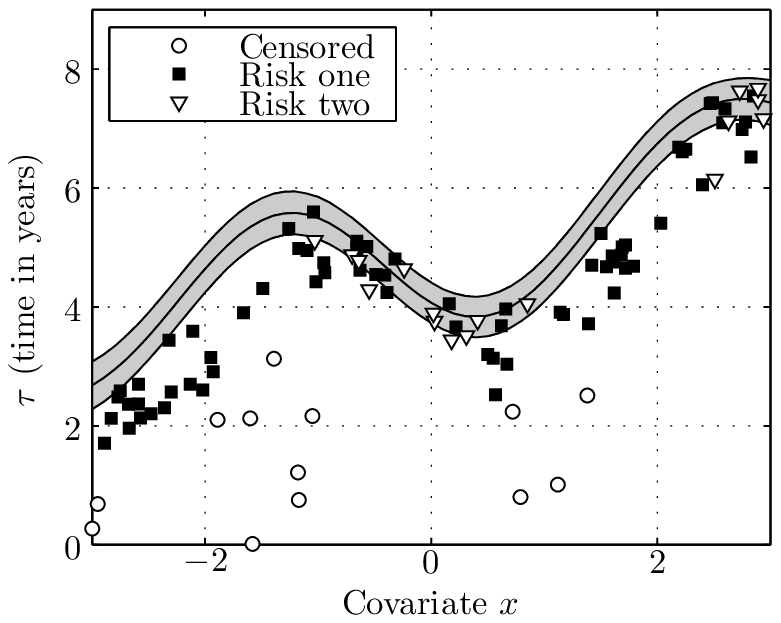}}\\
\subfloat[Independent risk 1]{\includegraphics[scale=0.7]{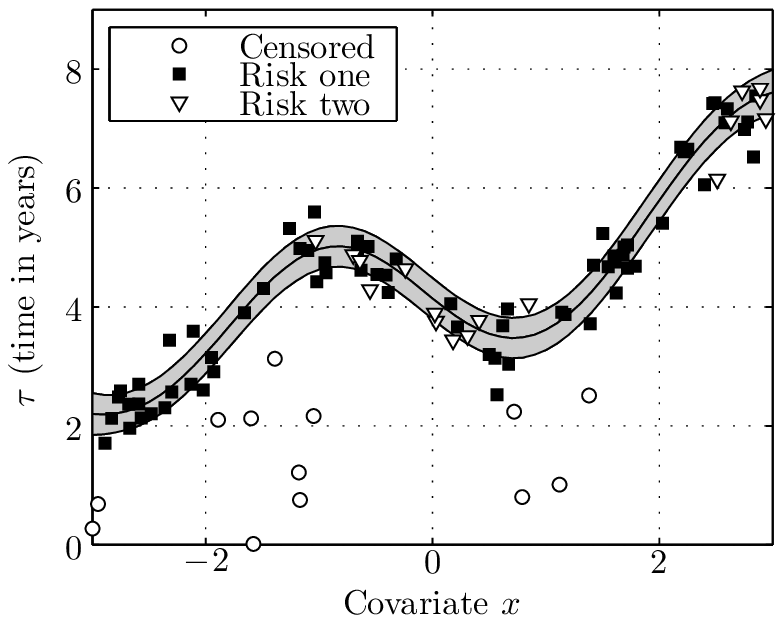}} &\subfloat[Independent risk 2]{\includegraphics[scale=0.7]{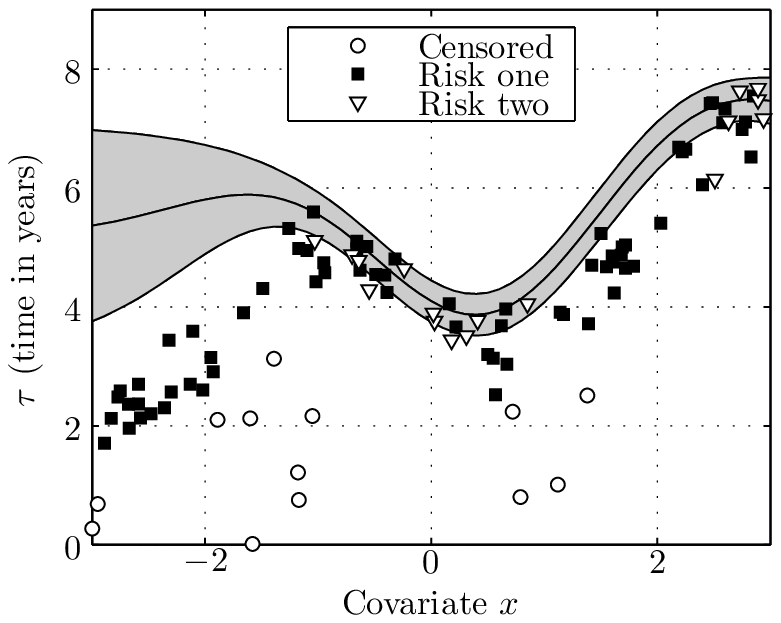}} \\
\end{tabular}
\caption[Comparison of dependent and independent GP regression on simulated data]{Results from simulated data with two competing risks. In (a) and (b) are inferred functions from both risks using a multiple output GP model with dependent risks allowed. Values of $\sigma = 0.2$ and $\omega=1.63$ were found which indicate strongly dependent risks. In (c) and (d) the inferred risks are forced to be independent by setting $\omega=0$. A value of $\sigma = 1.9$ is found which is sensible since the unique part of each risk must account for all of the variance. The advantage of assuming dependency can around $x\in(-3,-1)$ in (b) and (d). There are no risk two observations in this region and in the independent model there is much greater uncertainty since information form risk one cannot be utilised.}
\label{multi:fig:results3}
\end{figure}

In order to test the performance of GP regression in the presence of competing risks we generated survival data with dependent risks and compare two GP models, one which allows for dependency between risks and one with independent risks. Recall that by fixing the value of $\omega = 0$ we force the two risks to be independent. The results are shown in Figure \ref{multi:fig:results3}.

In Figure \ref{multi:fig:results3} (a) and Figure \ref{multi:fig:results3} (b) are results from a GP model where $\omega$ is inferred from the data. Hyperparameters were found to be $(\eta,\mu,\beta,\sigma,\omega,l) = (4.59, 0.41, 0.33, 0.20, 1.63, 1.01)$. The higher value of $\omega$ indicates that the model is assuming strong dependence between risks. In (c) and (d) are results from a second GP model with $\omega = 0$. Remaining hyperparameters were found to be $(\eta,\mu,\beta,\sigma,l) = (13.03, -0.60, 1.02, 1.90, 1.02)$. Note that the value of $\sigma$ is now higher as the unique part of each risk must explain all of the output variance. The advantage of allowing dependent risks becomes apparent when we examine the inferred risk 2 function towards the left of (b) and (d). In the independent model the uncertainty associated with the underlying function is much greater since knowledge of risk 1 is unavailable. In the dependent model a more accurate recovery of the risk 2 function is obtained and the uncertainty is smaller since information from risk 1 events can be utilised more effectively. In this dataset the generating hyperparameters were $(\eta,\mu,\beta,\sigma,\omega,l) = (5, 0.5, 0.5, 0.5, 2.0, 1.0)$ so the dependent model performs much better, particularly towards the left of the $x$-axis despite the complete lack of risk 2 observations. Of course, with real data we will not have the luxury of knowing whether an assumption of dependence is correct or not but this example nevertheless illustrates the potential usefulness of our approach.

\subsection{Application to multi-dimensional covariates}
\label{multi:sec:highdim}

It is straightforward to deal with more than one covariate. Using the ARD hyperparameters outlined in Section \ref{sec:single} we can determine which covariates have the biggest impact on survival outcomes. In the Supporting Information we give an example of competing risks data with two covariates. Inferred ARD hyperparameters are $(l_1, l_2) = (0.52, 1.47)$ indicating that the first covariate is more important.

\section{Discussion and conclusion}
\label{sec:disc}

We have pursued an alternative route to many existing survival analysis methods --- which assume a parametric or semi-parametric hazard rate --- and focus on directly characterising the relationship between covariates and event times in a flexible and non-parametric manner. GP regression provides a powerful and elegant means to achieve this. All relevant quantities are inferred in Bayesian manner and we can obtain probabilistic predictions, survival probabilities and hazard rates. It is straightforward to incorporate censored or truncated observations (or combination thereof). We found that the GP hazard rate model used by \cite{JOE12} also performed well on non-monotonic data but the hyperparameters are not easy to interpret. This is because the inferred function represents the relationship between the covariates and hazard rate whereas in our case the inferred function is conceptually more straightforward and easier to interpret.

An interesting question to consider is how to interpret the underlying function we infer. In standard GP regression the function values would be considered `noise-free' outputs which are then corrupted by Gaussian observational noise. The corresponding interpretation in our case would be that the functions represent a `noise-free' event time. However, Gaussian noise is not appropriate in that case since it is generally not plausible to claim events could be randomly reported before they actually occur. A more appropriate choice would be noise with a semi infinite support on $(0,\infty]$ that would represent a delay between the event occurring and the time it is diagnosed or recorded. In that case the underlying function could be interpreted as a noise-free event time. An alternative interpretation of the noise, however, is that it represents a disparity or mismatch between the assumed model and the actual model of the data. In that case it is acceptable to regard the function values as noise-free event times. In future work it may be interesting to explore alternative noise distributions.

Multiple output GP regression provided a natural route to incorporate competing risks data. Working with the event time density --- also called the latent event time approach --- has been criticised for a number of reasons. The most serious objection is that the joint event time density is unobservable and cannot be inferred from the observed failure and censoring times (due to the identifiability problem). While one may assume that the event times are independent or perhaps assume some parametric density to model dependencies one cannot conclude on the basis of the observed data alone whether or not the event times are independent. Some authors such as \cite[Section 3.3]{BEYER12} have claimed that the latent failure times lack plausibility or are too hypothetical in nature since we are positing the existence of quantities which in reality can never be measured (see \cite[Section 8.2]{KALB02} for further in depth discussion).

We argue instead that the latent failure time approach does provide a useful conceptual framework. There are two aspects in particular where such an approach is useful. The first is in making predictions for new patients (who are still alive) since the time until the different events is highly relevant. The second is in estimating what happens when one or more risks are disabled. In both cases the marginal survival functions are the relevant quantities. In our GP formulation both predictions and marginal survival probabilities are straightforward to compute due to the conditional independence of the event times. If we want to model dependent risks (despite not being able to test our assumption due to identifiability issues) then modelling the joint event time density is a convenient starting point. The fact that the data we observe in reality do not allow a direct view of this joint density does not mean that it is not a useful concept.

We have illustrated that GP regression provides a flexible method of analysing survival data. We impose few structural assumptions due to the non-parametric nature of GP regression and avoiding specification of the hazard rate. Future research could involve applying sparse GP regression techniques to our model in order to achieve greater computational efficiency.

\section*{Acknowledgments}
This work was funded under the European Commission FP7 Imagint Project, EC Grant Agreement number 259881.

\newpage
\section*{Supporting Information}
In these appendices we provide some additional background theory and results on GP regression for survival analysis. We outline our implementation of the GP hazard rate model that was used in \cite{JOE12}. Further results from simulated monotonic data (with a single risk and competing risks) are presented. Finally, partial derivates which are required for the Laplace approximations are given along with some implementational notes.
\appendix

\section{The GP hazard rate model}
\label{sec:gphazard}

For the purposes of comparison we wish to implement a GP hazard rate model recently used by \cite{JOE12} and \cite{SAV11}. In this section we outline the details of the model and how we implemented it.

\subsection{Model definition}

We compare our GP model to an alternative model that assumes a more traditional hazard rate $\pi_i(\tau|\vecx_i) = \lambda_0(\tau)\exp(f(\vecx_i))$ and a GP prior $p(\vecf|\matX)$ over the function values. Such models were used by \cite{JOE12} with a piecewise log-constant base hazard rate. The Laplace approximation was constructed in order to estimate the marginal likelihood and infer hyperparameters. The same model is discussed in \cite{SAV11}. We will use a baseline hazard rate corresponding to the Weibull distribution $\lambda_0(\tau) = \nu \tau^{\nu-1}$ with $\nu>0$, which implies $\Lambda_0(\tau) = \tau^{\nu}$. The likelihood contribution for individual $i$ is written in terms of the hazard rate $P(\tau_i,\Delta_i|f_i) = \pi_i(\tau)^{\Delta_i}e^{-\int_0^{\tau}\text{d}s\,\pi_i(s)}$. The posterior is obtained using Bayes' theorem:
\begin{equation}
p(\vecf|\matX,D,\tv)= \frac{p(D|\vecf,\tv)p(\vecf|\matX,\tv)}{\int\text{d}\vecf'\,p(D|\vecf',\tv)p(\vecf'|\matX,\tv)}.
\label{gpsa:eq:bayes}
\end{equation}
The negative log posterior is
\begin{align}
\mathcal{L}(\vecf) &= -\frac{1}{N}\sum_{i:\Delta_i=1}\bigg[\log\lambda_0(\tau_i) + f(\vecx_i)\bigg] +\frac{1}{N}\sum_{i=1}^N\Lambda_0(\tau_i)e^{f(\vecx_i)} \nonumber\\
&\qquad\qquad\qquad+\frac{1}{2N}\vecf\cdot\matK^{-1}\vecf+ \frac{1}{2N}\log|\matK| + \frac{1}{2}\log2\pi.
\end{align}
Similarly to the main text, hyperparameters are determined using the Laplace approximation of the marginal likelihood (partial derivatives are given bellow in Section \ref{app:gpsa:joensuu}). The predictive density over $f^*$ is Gaussian with mean and variance given by 
\begin{align}
\hat{\mu} & =\veck^*\cdot\matK^{-1}\hat{\vecf}\label{gpsa:eq:pred_mean}\\
\hat{\kappa} & = k(\vecx^*,\vecx^*) -\veck^*\cdot(\matK+\matW^{-1})^{-1}\veck^*.
\end{align}
The event time density is
\begin{equation}
p(\tau|f(\vecx_i)) = \lambda_0(\tau)e^{f(\vecx_i)}e^{-\Lambda_0(\tau)e^{f(\vecx_i)}}
\label{gpsa:eq:pred_etd}
\end{equation}
which can be used to obtain the predictive distribution over the event time:
\begin{equation}
p(\tau^*|\vecx^*,\matX,D) = \int \text{d}f^*\,\lambda_0(\tau^*)e^{f^*}e^{-\Lambda_0(\tau^*)e^{f^*}}\frac{e^{-\frac{1}{2\hat{\kappa}^2}(f^*-\hat{\mu})^2}}{(2\pi\hat{\kappa}^2)^{1/2}}.
\label{gpsa:eq:pred_time}
\end{equation}
This expression was found to be problematic in that $\int\text{d}\tau^*\,p(\tau^*|\vecx^*,\matX,D)\neq 1$. This is an example of non-commuting limits since if we first integrate the integrand in (\ref{gpsa:eq:pred_time}) with respect to $\tau^*$ and then $f^*$ we obtain a value of one. A rough explanation of why this occurs is that negative values of $f^*$ correspond to a protective effect, that is, the event time density places more probability mass away from the origin. As the value of $f^*$ decreases more probability mass is placed further and further away from the origin. In the limit $f^*\to-\infty$ then $\text{prob}(\tau^* =\infty) \to 1$. We have not produced a more formal examination of this issue but experience suggests that numerical computation of the mean and variance of (\ref{gpsa:eq:pred_time}) is infeasible. Consequently, the predictive mean $\left<\tau^*\right>$ and variance $\left<(\tau^*)^2\right> - \left<\tau^*\right>^2$ will be computed numerically from (\ref{gpsa:eq:pred_etd}). Note that this will underestimate the uncertainty since we are not taking into account the uncertainty in $f^*$.

\subsection{Application to interval censored data}
\label{gpsa:sec:joe12interval}

The GP hazard rate model can also accommodate interval censored observations. We can use
\begin{align}
\mathcal{L}(\vecf) &= -\frac{1}{N}\sum_{i:\Delta_i=1}\log[S(t_i^l|f_i) - S(t_i^u|f_i)] -\frac{1}{N}\sum_{i:\Delta_i=0}\log S(t_i|f_i) -\frac{1}{N}\log p(\vecf|\matX,\tv)
\end{align}
 with $S(\tau|f_i) = \exp(-\Lambda_0(\tau)e^f_i)$ to determine the optimal values of $\vecf$. Hyperparameters are found using the Laplace approximation.

\section{Additional results and figures}

A plot of the event time transformation is given. We then provide an example of our GP regression method on monotonic survival data and compare the survival curves and hazard rates to those obtained with a WPHM or Cox model. We apply the GP hazard rate model to the single risk simulated data given in the main text. We then apply our multiple output GP regression to monotonic competing risks data and compare its performance to the WPHM. Finally, we give an example of using the ARD hyperparameters with two dimensional covariates and competing risks.

\subsection{The time-to-event transformation}

\begin{figure}[htb!]
\centering
\includegraphics[scale = 0.8]{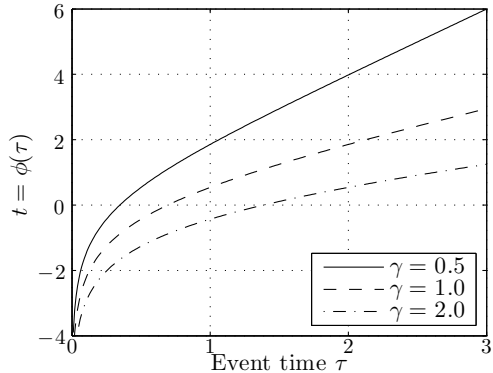}
\caption[Time-to-event transformation in GP regression]{Plot of the time transformation $\phi(\tau)$ that was used in the GP implementation. Note that for values of $\tau>\gamma$ the transformation is approximately linear. By adjusting the value of $\gamma$ such that it is less than the earliest observed event time we effectively end up with a linear mapping. The effect of the transformation can be seen when negative values of $t = \phi(\tau)$ are predicted since they are `squashed' into the positive half of the real line.}
\label{gpsa:fig:transform}
\end{figure}

In Figure \ref{gpsa:fig:transform} is a plot of the event time transformation $\phi(\tau) = \log(e^{\tau/\gamma}-1)$.

\begin{figure}[hbt!]
\centering
\begin{tabular}{c c}
\subfloat[WPHM]{\includegraphics[scale=0.8]{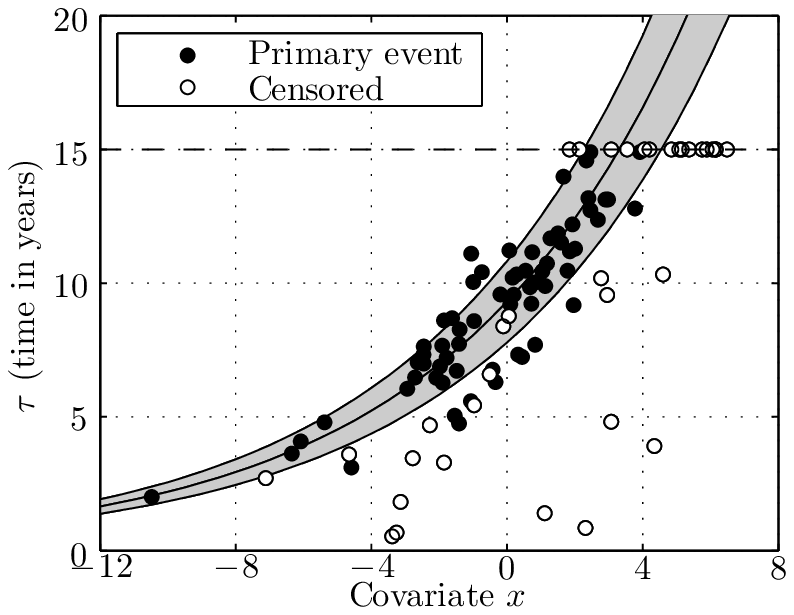}} & \subfloat[Our model]{\includegraphics[scale = 0.8]{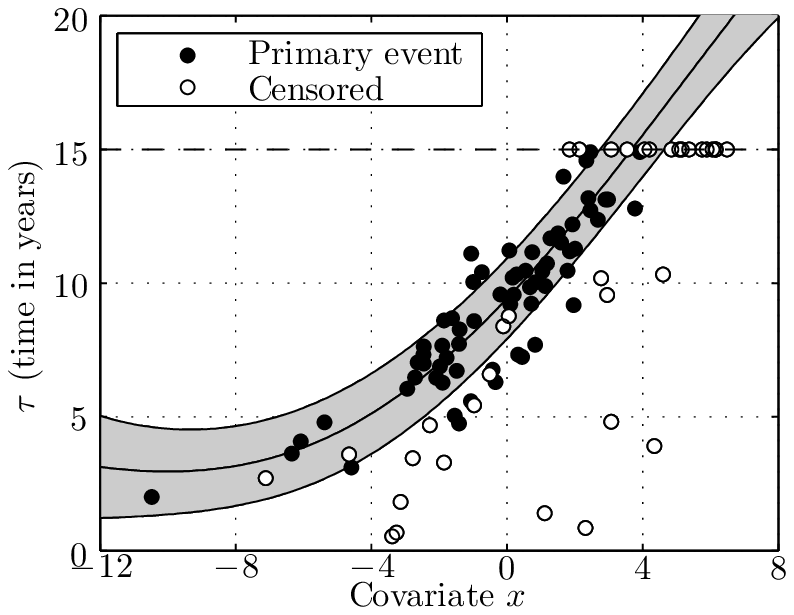}}\\
\subfloat[Survival curves]{\includegraphics[scale=0.8]{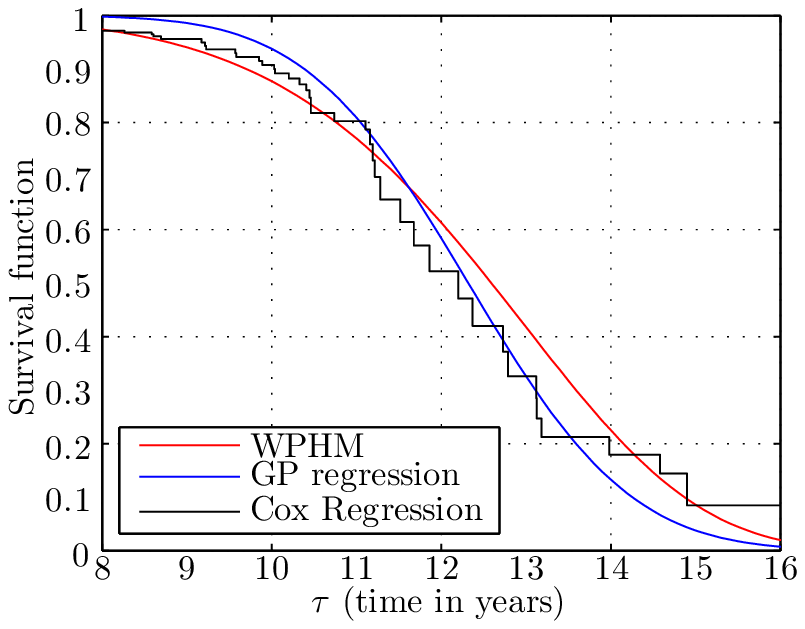}} &\subfloat[Hazard rates]{\includegraphics[scale=0.8]{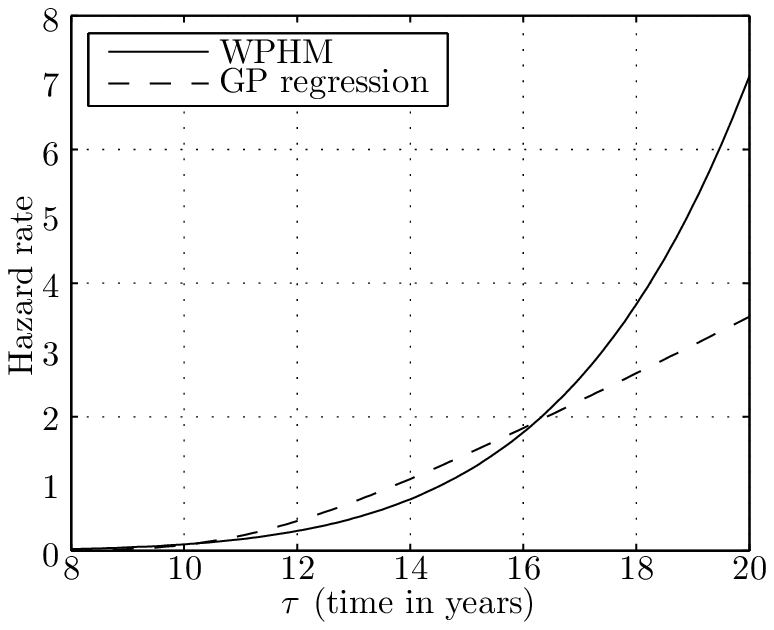}} \\
\end{tabular}
\caption[Comparison of GP regression to the WPHM on monotonic simulated data]{Example of data generated according to the WPHM assumptions. In (a) are results from fitting the WPHM. We found $(\beta,\rho,\nu) = (-1.04,9.93, 7.23)$. In (b) are results from our GP model which is also capable of handling these monotonic data although the uncertainty is greater towards the left of the figure. In (c) we compare survival curves corresponding to an individual with $x=2$ from the WPHM, our GP model and a Cox proportional hazards model (with $\beta_{cox} = -1.02$). In (d) are hazard rates  from the WPHM and our GP model for an individual with $x=2$.}
\label{gpsa:fig:example3}
\end{figure}

\subsection{Simulated monotonic data with a single risk}

Here we generated simulated data corresponding to the WPHM. We begin by choosing values of $\bv$, $\rho$, $\nu$ manually. Covariate vectors $\vecx_i$ are generated from a uniform distribution on a finite region of the covariate space. Event times are generated using the inverse of the cumulative distribution:
\begin{equation}
C_i(\tau) = 1-e^{-\Lambda_0(\tau)e^{\bv\cdot\vecx_i}}.
\end{equation}
Random numbers $z\in[0,1]$ are generated from a uniform density. An event time corresponding to $\vecx_i$ is
\begin{equation}
\tau_i = \rho\left(-e^{-\bv\cdot\vecx_i}\log(1-z)\right)^{1/\nu}.
\label{gpsa:eq:synthgen}
\end{equation}
Finally independent censoring is simulated by randomly selecting a subset of the individuals and generating a random number from a uniform distribution defined on the interval $[0,\tau_i)$ which is then recorded as the time of censoring.

These data are shown in Figure \ref{gpsa:fig:example3} and have a monotonic relationship between the event time and the covariate. We ran both the WPHM and our GP model in order to see how our model performs on data that can readily be analysed with existing tools. In Figure \ref{gpsa:fig:example3} (a) are the results from running the WPHM. Visually, it is clear that the model achieves a good fit. In Figure \ref{gpsa:fig:example3} (b) are the results from our GP model. Our model has also achieved a good fit. One difference between both models is that the GP model has greater uncertainty towards the left of the figure. This is appropriate since there are very few observations here so consequently our knowledge of the underlying function is less firm.

In Figure \ref{gpsa:fig:example3} (c) we have compared the survival functions of the WPHM, our GP model, and a Cox proportional hazards model. The survival functions all correspond to an individual with $x=2$. It is clear that all three models are giving broadly similar survival probabilities. Finally, in Figure \ref{gpsa:fig:example3} (d) we plot the hazard rates corresponding to an individual with $x=2$ for both the WPHM and our GP model. Note that the hazard rate in our GP model is approximately linear for large times.

\begin{figure}[hbt!]
\centering
\begin{tabular}{c c}
\subfloat[Right censoring only]{\includegraphics[scale=0.8]{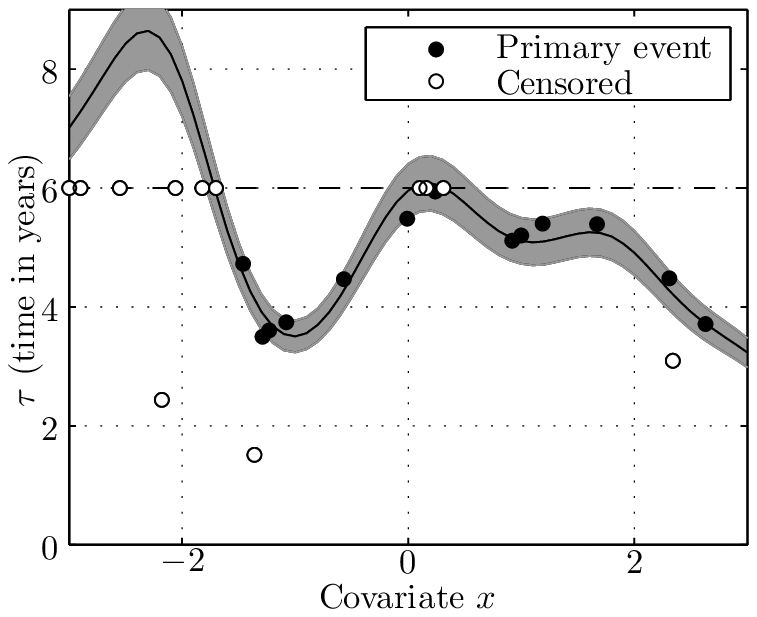}} &\subfloat[Interval censoring]{\includegraphics[scale=0.8]{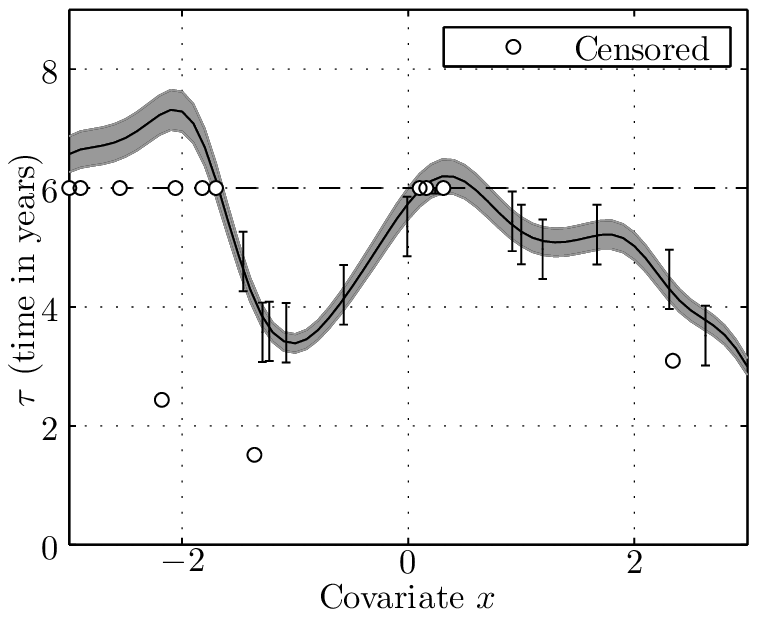}} \\
\end{tabular}
\caption{In (a) is the inferred function from data in Figure 1 of the main text using the GP hazard rate model. The inferred hyperparameters are $(\eta,\beta,\sigma,l) = (-45.5,27.2,64.7,0.47)$ although these are difficult to interpret since the underlying function appears as $\exp(f(x))$ in the hazard rate. In (b) are results from the interval censored version of these data. Note that the uncertainty is underestimated in this model.}
\label{gpsa:fig:example2}
\end{figure}

\subsection{Application of the GP hazard rate model to right censored and interval censored simulated data}

Here we apply the GP hazard rate model described above in Section \ref{sec:gphazard} to the simulated data presented in Figure 1 of the main text. The results are shown in Figure \ref{gpsa:fig:example2}. Note that the hyperparameters are more difficult to interpret and that the uncertainty is underestimated.

\subsection{Monotonic competing risks survival data with dependent competing risks}

\begin{figure}[h!]
\centering
\begin{tabular}{c c}
\subfloat[GP model inferred risk 1]{\includegraphics[scale=0.8]{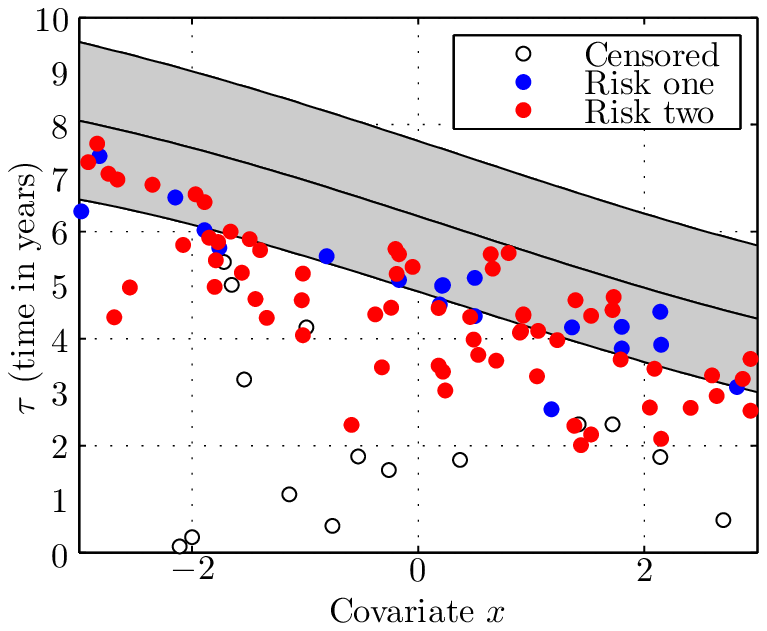}} & \subfloat[GP model inferred risk 2]{\includegraphics[scale = 0.8]{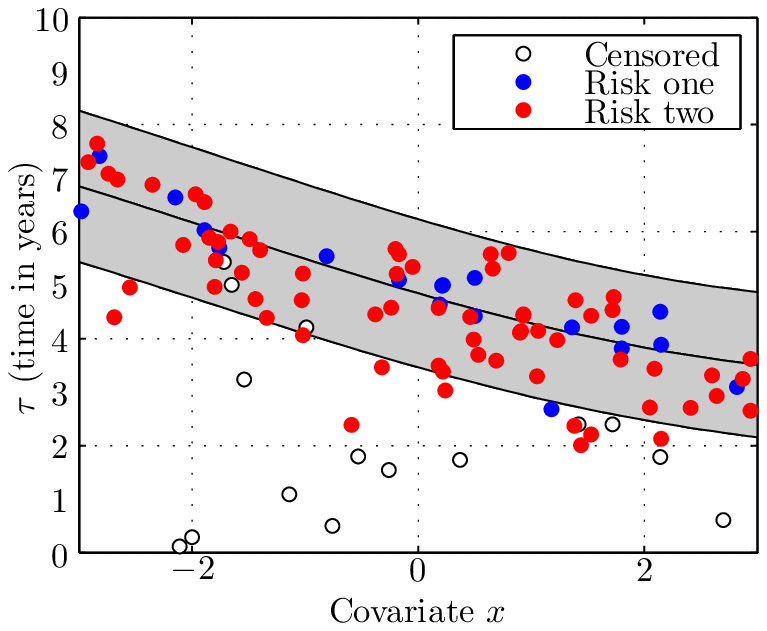}}\\
\subfloat[WPHM inferred risk 1]{\includegraphics[scale=0.8]{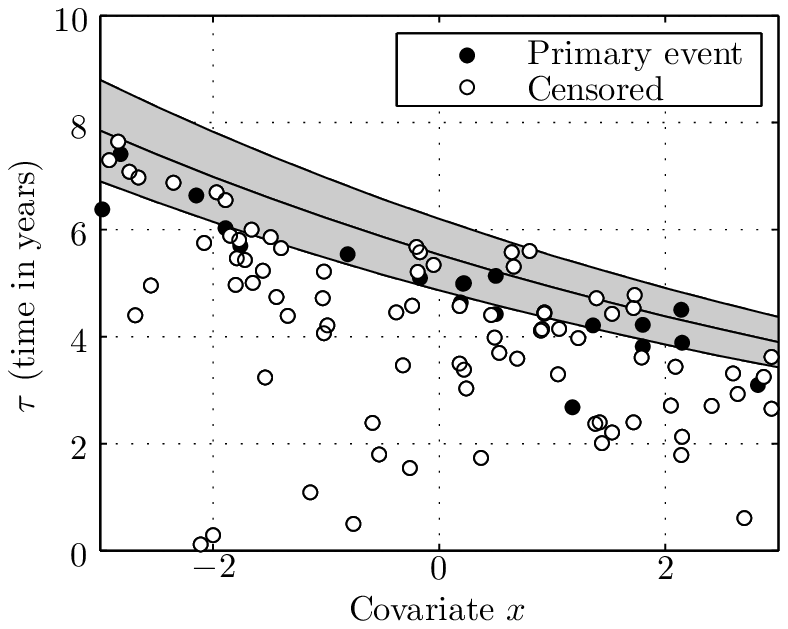}} &\subfloat[WPHM inferred risk 2]{\includegraphics[scale=0.8]{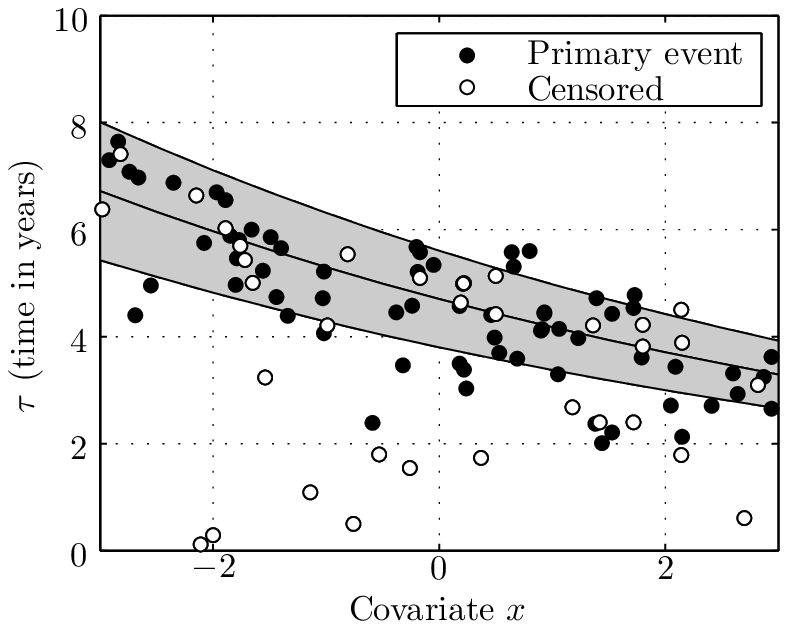}} \\
\end{tabular}
\caption[Example of GP regression on simulated monotonic competing risks data]{Example of monotonic survival data with two dependent competing risks. In Figures (a) and (b) are plots of the inferred risk one and risk two functions respectively using the GP model. Inferred hyperparameters are $(\eta,\mu,\beta,\sigma,\omega,l) = (5.95,2.16,1.40,0,0.96,4.56)$. The value $\omega\neq 0$ reflects the fact that the model has assumed the risks are dependent. The characteristic length scale $l=4.56$ indicates that the function is changing relatively slowly with respect to the covariate. In Figures (c) and (d) are results from running two independent WPHM models. In each model only one of the risks is regarded as the primary risk and all other events are considered as right censored.}
\label{multi:fig:results2}
\end{figure}

Here we study the performance of multiple output GP regression on a univariate dataset with monotonic competing risks event times. We train a GP model and compare this to two independent WPHM models (Figure \ref{multi:fig:results2}). In this case values of $\omega = 0.95$ and $\sigma = 0$ were inferred which indicates that the risks are completely dependent with no variance due to unique components. The characteristic length scale $l=4.56$ illustrates that the function is changing relatively slowly with respect to the covariate. In Figure \ref{multi:fig:results2} (a) we see that the GP model has inferred a risk 1 function that lies `above' most of the observed event times. This is not unreasonable since all of the risk 2 events are effectively censoring events from the point of view of risk 1. Therefore we know that risk 1 events must occur after risk 2 event times. Consequently the data likelihood is maximised by placing the risk 1 function slightly `above' the risk 2 events. A similar effect can be seen in Figure \ref{multi:fig:results2} (c) since the risk 2 events are also regarded as censoring events in the WPHM. We compute the MSE in a validation set of 100 samples. Results are shown in Table \ref{multi:tab:results2}. The GP model performs slightly worse that the WPHM model, particularly when it comes to predicting risk 1 events. This appears to be due to the fact that the GP model has inferred a risk 1 function that is slightly `higher' than the WPHM risk 1 function. Consequently, the predicted risk 1 events are overestimating the time to event. In this dataset there are 66 risk 2 events compared to 19 risk 1 events which helps to explain the poorer predictions for risk 1. Also, because of the way the validation data are generated only the earliest risk 1 events are reported which leads to larger MSE values in both models because the predicted event time are tend to be later than the reported event times.

\begin{table}[h]
  \begin{center}
    \begin{tabular}{|c|c|c|}
    \hline
     & WPHM & GP regression \\
	\hline
	Risk 1 MSE (years$^2$)& 1.42 & 3.10\\
    Risk 2 MSE (years$^2$)& 0.75 & 0.84\\
    \hline
    \end{tabular}
  \end{center}
\caption[GP regression applied to monotonic simulated competing risks data]{Comparison of mean square error (MSE) between the WPHM and the GP model on 100 validation samples corresponding to the training data in Figure \ref{multi:fig:results2}. The GP model has slightly poorer performance than the WPHM, particularly on risk 1. See the main text for further discussion.}
\label{multi:tab:results2}
\end{table}

\subsection{Example of two dimensional covariates}

\begin{figure}[h!]
\centering
\includegraphics[scale=0.6]{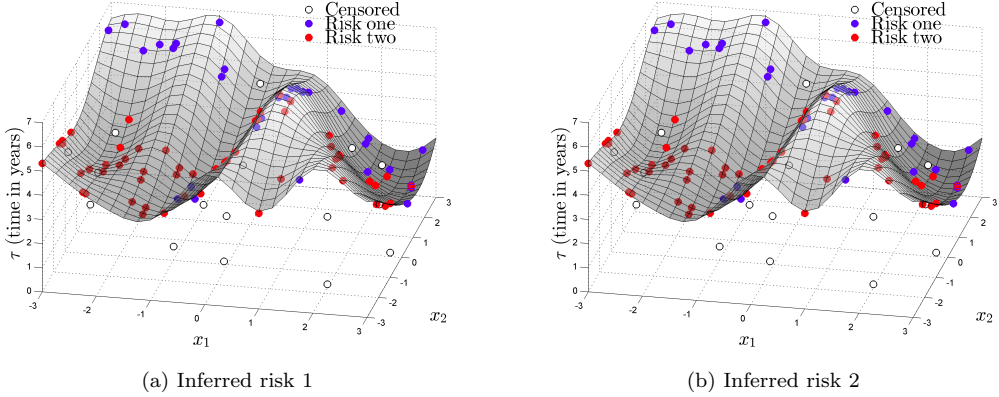}
\caption[GP regression with two dimensional covariates]{Example of GP regression with two dimensional covariates and strongly dependent risks. Covariate $x_1$ has an inferred characteristic length scale of $l_1=0.52$ compared to $l_2=1.47$. The values used to generate the data were $(l_1,l_2) = (0.5, 1.5)$. This indicates that the first covariate is more relevant to determining survival outcomes. This is reflected in the plots since the function is more variable in the $x_1$ direction. The remaining hyperparameters were found to be $(\eta,\beta,\sigma,\omega) = (4.83, 0.17, 0.23, 0.89)$. Note that higher value of $\omega$ reflects the fact that the GP model has assumed a strong dependence between risks.}
\label{multi:fig:results4}
\end{figure}

The model is fully capable of dealing with multi-dimensional covariates. By using the squared exponential kernel with automatic relevance determination (ARD) hyperparameters we can determine which covariates are the most important. This is analogous to examining the regression coefficients in a Cox model to see which covariates have the greatest impact on survival outcomes. In the example shown in Figure \ref{multi:fig:results4} we find that $(l_1, l_2) = (0.52, 1.47)$ indicating that the first covariate is more important.

\section{Numerical implementation}
\label{app:sec:gplvm}

\subsection{GP regression with a single risk and independent censoring}
\label{gpaft_num}

A number of numerical issues arise during the implementation our GP model. Numerical instability can occur when computing the negative log likelihood function (with right censoring). The problematic terms are the hazard rates $\pi(t|f) = p(t|f)/S(t|f)$ where $p(t|f)$ is a Gaussian density and $S(t|f)$ is the corresponding survival function. The hazard rates do not appear in negative log likelihood but the same terms do occur in the partial derivatives (see (\ref{app:eq:prob1}) and (\ref{nablarho})). This quantity can be written in terms of the \emph{complementary error function} 
\begin{equation}
\text{erfc}(x) = \frac{2}{\sqrt{\pi}}\int_x^{\infty}\text{d}s\,e^{-s^2}.
\end{equation}
If we define $h = (t-f)/\beta\sqrt{2}$ then
\begin{align}
\log S(t|f) &= \log\left(\tfrac{1}{2}\text{erfc}(h)\right)\label{gpsa:eq:erfclog}\\
\frac{\partial}{\partial f}\log S(t|f) &= \frac{2}{\sqrt{2\pi}\beta}\frac{e^{-h^2}}{\text{erfc}(h)}\label{gpsa:eq:erfcd}\\
\frac{\partial^2}{\partial f^2}\log S(t|f) &= -\left(\frac{2}{\sqrt{2\pi}\beta}\frac{e^{-h^2}}{\text{erfc}(h)}\right)^2 +\frac{\sqrt{2}h}{\beta}\left(\frac{2}{\sqrt{2\pi}\beta}\frac{e^{-h^2}}{\text{erfc}(h)}\right).\label{gpsa:eq:erfcdd}
\end{align}
The hazard rate is also given by (\ref{gpsa:eq:erfcd}). For large $h$ the quantity $e^{-h^2}/\erfc(h)$ becomes numerically unstable since both numerator and denominator tend towards zero. This is solved by using the asymptotic expansion of the complementary error function \cite{MEN60}:
\begin{equation}
\text{erfc}(h) = \frac{e^{-h^2}}{h\sqrt{\pi}}\left[1 - \frac{1}{2h^2} + \frac{2}{(2h^2)^2} - \frac{8}{(2h^2)^3} + \cdots\right]\quad\text{for $h \gg 0$.}
\end{equation}
This results in the following approximations which are numerically stable:\footnote{Using the approximations for $h>20$ gives acceptable performance in Matlab.}
\begin{align}
\log S(t|f) &= -h^2 -\log h -\log 2\sqrt{\pi}+\log\left[1 - \frac{1}{2h^2} + \frac{2}{(2h^2)^2} - \frac{8}{(2h^2)^3} + \cdots\right]\\
\frac{\partial}{\partial f}\log S(t|f)&= \frac{h\sqrt{2}}{\beta} \left[1 - \frac{1}{2h^2} + \frac{2}{(2h^2)^2} - \frac{8}{(2h^2)^3} + \cdots\right]^{-1}\label{gpsa:eq:approxerfcd}.
\end{align}
The approximation (\ref{gpsa:eq:erfcd}) can be substituted directly into (\ref{gpsa:eq:erfcdd}) to obtain an approximation for the second order partial derivative. Note that the hazard rate is approximately linear for large $h$ which is consistent with Figure \ref{gpsa:fig:example3} (d). 

\subsection{The GP hazard rate model for a single risk with independent censoring}
\label{gphr_num}

Some numerical issues arise in the computation of $\log(S(\tau_i^l)-S(\tau_i^u))$ while we search the parameter space for an optimal solution (and when we compute the partial derivatives (\ref{app:eq:problematique1}) and (\ref{app:eq:problematique2})) using the GP hazard rate model. This is because $S(\tau) =\exp(-\Lambda_0(\tau)e^f)$ can take extremely small values and unlike the right censored case the $\log$ does not cancel the exponentials because of the sum. In this case we write
\begin{align}
\log(e^{-x_1} - e^{-x_2}) & = \log\{e^{-x_1}(1-e^{x_1-x_2})\}\\
& = \left\{\begin{array}{ll}
 -x_1 +\log(1-e^{x_1-x_2})&\qquad\text{when $-C\leq x_1-x_2$}\\
 -x_1 -e^{x_1-x_2} -\tfrac{1}{2}e^{2(x_1-x_2)}&\qquad\text{when $x_1-x_2<-C$.}\\
 \end{array}
 \right.
\end{align}
Note that $x_2\geq x_1$ in the context of this model. The constant $C$ is a cutoff that depends on the numerical accuracy of implementation.\footnote{In Matlab $C=10$ was found to be sufficient.} Furthermore, a similar problem occurs with the computation of first and second order gradients. A similar trick can rectify the problem. The offending terms from the gradient in Section \ref{app:gpsa:interval_joensuu} are (\ref{app:eq:problematique1}) and (\ref{app:eq:problematique2}) and can be rewritten as
\begin{align}
\frac{-x_1e^{-x_1}+x_2e^{-x_2}}{e^{-x_1} -e^{-x_2}} = \frac{-x_1+x_2e^{x_1-x_2}}{1-e^{x_1-x_2}},
\end{align}
and the second order derivatives are given by
\begin{equation}
-\delta_{ij}\left(\frac{-x_1+x_2e^{x_1-x_2}}{1-e^{x_1-x_2}}\right)^2 +\delta_{ij}\frac{-x_1+x_1^2+(x_2-x_2^2)e^{x_1-x_2}}{1-e^{x_1-x_2}}.
\end{equation}

\section{Partial derivatives}

In this section we derive the first and second order partial derivatives of the likelihood functions corresponding to the various GP regression methods used in the main text and here. The partial derivatives are required for gradient based optimisation of the likelihood and construction of the Laplace approximation of the marginal likelihood. For clarity we will rewrite the corresponding negative log likelihood from each section.

\subsection{GP regression with a single risk}
\label{app:gpsa:partials}

The negative log likelihood function is
\begin{equation}
\mathcal{L}(\vecf) =-\frac{1}{N}\sum_{i:\Delta_i=1}\log p(t_i|f_i) -\frac{1}{N}\sum_{i:\Delta_i=0}\log S(t_i|f_i) + \frac{1}{2N}(\vecf-\ev)\cdot\matK^{-1}(\vecf-\ev) + \frac{1}{2N}\log|\matK|.
\end{equation}
First order partial derivatives are
\begin{equation}
\frac{\partial}{\partial f_i}\mathcal{L}(\vecf) = -\frac{1}{N}\sum_{k:\Delta_k=1}\frac{\partial}{\partial f_i}\log p(t_k|f_k) -\frac{1}{N}\sum_{k:\Delta_k=0}\frac{\partial}{\partial f_i}\log S(t_k|f_k)+\frac{1}{N}(\matK^{-1}(\vecf-\ev))_i\label{nablapsi}
\end{equation}
where
\begin{equation}
\frac{\partial}{\partial f_i}\log p(t_k|f_k) = \delta_{\Delta_i,1}\delta_{ik}\beta^{-2}(t_k-f_k)
\end{equation}
and
\begin{align}
\frac{\partial}{\partial f_i}\log S(t_k|f_k) & =\delta_{\Delta_i,0}\delta_{ik}\frac{1}{S(t_k|f_k)}\int_{t_k}^{\infty}\text{d}s\,\frac{\partial}{\partial f_i}\frac{e^{-\frac{1}{2\beta^2}(s - f_k)^2}}{\sqrt{2\pi}\beta}\nonumber\\
& =\delta_{\Delta_i,0}\delta_{ik}\frac{1}{S(t_k|f_k)}\frac{1}{\sqrt{2\pi}\beta^3}\int_{t_k}^{\infty}\text{d}s\,(s-f_k)e^{-\frac{1}{2\beta^2}(s - f_k)^2}\nonumber\\
& = \delta_{\Delta_i,0}\delta_{ik}\frac{1}{S(t_k|f_k)}\frac{1}{\sqrt{2\pi}\beta} e^{-\frac{1}{2\beta^2}(t_k-f_k)^2}.
\label{app:eq:prob1}
\end{align}
Second order partial derivatives are
\begin{align}
\frac{\partial^2}{\partial f_i\partial f_j}\mathcal{L}(\vecf) &= -\frac{1}{N}\delta_{ij}\sum_{k:\Delta_k=0}\frac{\partial^2}{\partial f_i^2}\log p(t_k|f_k) -\frac{1}{N}\delta_{ij}\sum_{k:\Delta_k=1}\frac{\partial^2}{\partial f_i^2}\log S(t_k|f_k)+\frac{1}{N}\matK^{-1}_{ij}\nonumber\\
&=\frac{1}{N}(\matW + \matK^{-1})_{ij}
\end{align}
where the diagonal matrix $\matW$ is defined by $\matW_{ii} = -\frac{\partial^2}{\partial f_i^2}\log p(D|\vecf)$ with
\begin{equation}
\frac{\partial^2}{\partial f_i^2} \log p(t_k|f_k) = -\delta_{\Delta_i,1}\delta_{ik}\beta^{-2}
\end{equation}
and
\begin{align}
\frac{\partial^2}{\partial f_i^2}\log S(t_k|f_k) &= \delta_{\Delta_i,0}\delta_{ik}\left[-\left(\frac{1}{S(t_k|f_k)}\frac{1}{\sqrt{2\pi}\beta} e^{-\frac{1}{2\beta^2}(t_k-f_k)^2}\right)^2\right.\nonumber\\
&\qquad\qquad\qquad\qquad\left. + \frac{(t_k-f_k)}{\beta^2}\left(\frac{1}{S(t_k|f_k)}\frac{1}{\sqrt{2\pi}\beta} e^{-\frac{1}{2\beta^2}(t_k-f_k)^2}\right)\right]\label{nablarho}.
\end{align}
Note that the elements of $\matW$ are non negative because $S(t|f)$ is log concave. This implies that the Hessian is positive definite which we expect at the minimum of $\mathcal{L}(\vecf)$. To see that $S(t|f)$ is log concave we note that
\begin{align*}
S(t|f) & = \int_{t}^{\infty}\text{d}s\,\frac{e^{-\frac{1}{2\beta^2}(s - f)^2}}{\sqrt{2\pi}\beta} =  \int_{t-f}^{\infty}\text{d}s\,\frac{e^{-\frac{1}{2\beta^2}s^2}}{\sqrt{2\pi}\beta} =  \int_{-\infty}^{f-t}\text{d}s\,\frac{e^{-\frac{1}{2\beta^2}s^2}}{\sqrt{2\pi}\beta} 
\end{align*}
which is the cumulative distribution function for a Gaussian which is log concave.

\subsection{GP regression with interval censored data}
\label{app:gpsa:interval_partials}
The negative log likelihood function is
\begin{equation}
\mathcal{L}(\vecf) = -\frac{1}{N}\sum_{i:\Delta_i=1}\log [S(t_i^l|f_i) - S(t_i^u|f_i)] -\frac{1}{N}\sum_{i:\Delta_i=0}\log S(t_i|f_i) -\frac{1}{N}\log p(\vecf|\matX).
\end{equation}
For compactness we define the interval $I_k = (t^l,t^u)$ and
\begin{equation}
\psi(I_k|f_k) = S(\tau_i^l|f_k) - S(\tau_i^u|f_k) = \int_{t^l}^{t^u}\text{d}s\,\frac{e^{-\frac{1}{2\beta^2}(s-f_k)^2}}{\sqrt{2\pi}\beta}.
\end{equation}
The first order partial derivatives are
\begin{align}
\frac{\partial}{\partial f_i}\log \psi(I_k|f_k) &= \frac{1}{\psi(I_k|f_k)}\frac{\partial}{\partial f_i}\int_{t^l}^{t^u}\text{d}s\,\frac{e^{-\frac{1}{2\beta^2}(s-f_k)^2}}{\sqrt{2\pi}\beta}\nonumber\\
&=\frac{1}{\psi(I_k|f_k)}\frac{\partial}{\partial f_i}\int_{f_k-t^u}^{f_k-t^l}\text{d}s\,\frac{e^{-\frac{1}{2\beta^2}s^2}}{\sqrt{2\pi}\beta}\nonumber\\
&=\frac{1}{\psi(I_k|f_k)}\left(\frac{e^{-\frac{1}{2\beta^2}(t^l-f_k)^2}}{\sqrt{2\pi}\beta} - \frac{e^{-\frac{1}{2\beta^2}(t^u-f_k)^2}}{\sqrt{2\pi}\beta}\right).
\end{align}
The second order partial derivatives are
\begin{align}
\frac{\partial^2}{\partial f_i^2}\log \psi(I_k|f_k) &=-\left[\frac{1}{\psi(I_k|f_k)}\left(\frac{e^{-\frac{1}{2\beta^2}(t^l-f_k)^2}}{\sqrt{2\pi}\beta} - \frac{e^{-\frac{1}{2\beta^2}(t^u-f_k)^2}}{\sqrt{2\pi}\beta}\right)\right]^2\nonumber\\
& \qquad+ \frac{1}{\psi(I_k|f_k)}\left(\frac{(t^l-f_k)}{\beta^2}\frac{e^{-\frac{1}{2\beta^2}(t^l-f_k)^2}}{\sqrt{2\pi}\beta} - \frac{(t^u-f_k)}{\beta^2}\frac{e^{-\frac{1}{2\beta^2}(t^u-f_k)^2}}{\sqrt{2\pi}\beta}\right).
\end{align}
The partial derivatives ${\partial^2}/{\partial f_i \partial f_j} \log\psi =0$ for $i\neq j$. If we define $h^u = (t^u-f)/\beta\sqrt{2}$ and $h^l = (t^l-f)/\beta\sqrt{2}$ then
\begin{align}
\psi(I_k|f_k) &= \frac{1}{2}\erfc(h^l) - \frac{1}{2}\erfc(h^u)\\
\frac{\partial}{\partial f_i}\log \psi(I_k|f_k) &= \frac{2}{\sqrt{2\pi}\beta}\frac{e^{-(h^l)^2}-e^{-(h^u)^2}}{\erfc(h^l)-\erfc(h^u)}\\
\frac{\partial^2}{\partial f_i^2}\log \psi(I_k|f_k) &= -\left(\frac{2}{\sqrt{2\pi}\beta}\frac{e^{-(h^l)^2}-e^{-(h^u)^2}}{\erfc(h^l)-\erfc(h^u)}\right)^2 + \frac{2}{\sqrt{\pi}\beta}\frac{h^le^{-(h^l)^2} - h^ue^{-(h^u)^2}}{\erfc(h^l)-\erfc(h^u)},
\end{align}
and we can use the asymptotic expansion of the complementary error function to avoid any numerical difficulties (see Section \ref{gpaft_num}).
\subsection{The GP hazard rate model}
\label{app:gpsa:joensuu}
The negative log likelihood is
\begin{align}
\mathcal{L}(\vecf) &= -\frac{1}{N}\sum_{i:\Delta_i=1}\bigg[\log\lambda_0(\tau_i) + f(\vecx_i)\bigg] +\frac{1}{N}\sum_{i=1}^N\Lambda_0(\tau_i)e^{f(\vecx_i)} +\frac{1}{2N}\vecf\cdot\matK^{-1}\vecf\nonumber\\
&\qquad\qquad\qquad\qquad\qquad + \frac{1}{2N}\log|\matK| + \frac{1}{2}\log2\pi.
\end{align}
First order derivatives are
\begin{equation}
\frac{\partial}{\partial f_i}\mathcal{L}(\vecf) = -\frac{1}{N}\delta_{1,\Delta_i} + \frac{1}{N}\Lambda_0(\tau_i)e^{f(\vecx_i)} + \frac{1}{N}(\matK^{-1}\vecf)_i.
\end{equation}
Second order partial derivatives are
\begin{align}
\frac{\partial^2}{\partial f_i\partial f_j}\mathcal{L}(\vecf) &= \frac{1}{N}\delta_{ij}\Lambda_0(\tau_i)e^{f(\vecx_i)}+\frac{1}{N}\matK^{-1}_{ij}\nonumber\\
&=\frac{1}{N}(\matW + \matK^{-1})_{ij}
\end{align}
where $\matW$ is a diagonal matrix defined by $\matW_{ii} = \Lambda_0(t_i)e^{f(\vecx_i)}$.

\subsection{The GP hazard rate model with interval censoring}
\label{app:gpsa:interval_joensuu}
The negative log likelihood function is
\begin{align}
\mathcal{L}(\vecf) &= -\frac{1}{N}\sum_{i:\Delta_i=1}\log S(\tau_i|f_i)-\frac{1}{N}\sum_{i:\Delta_i=0}^N\log [S(\tau_i^l|f_i)-S(\tau_i^u|f_i)] +\frac{1}{2N}\vecf\cdot\matK^{-1}\vecf\nonumber\\
&\qquad\qquad\qquad\qquad\qquad + \frac{1}{2N}\log|\matK| + \frac{1}{2}\log2\pi.
\end{align}
The first order partial derivatives can be obtained from
\begin{equation}
-\frac{1}{N}\sum_{k:\Delta_k=1}\frac{\partial}{\partial f_i}\big[-\Lambda_0(\tau_k)e^{-f(\vecx_k)}\big] = \frac{1}{N}\Lambda_0(\tau_i)e^{f(\vecx_i)}
\end{equation}
and
\begin{align}
\sum_{k:\Delta_k=0}\frac{\partial}{\partial f_i}\log [S(\tau_k^l)-S(\tau_k^u)] &= \frac{1}{S(\tau_i^l) - S(\tau_i^u)}\left(-\Lambda_0(\tau_i^l)e^{f(\vecx_i)}e^{-\Lambda_0(\tau_i^l)e^{f(\vecx_i)}}\right.\nonumber\\
&\qquad\qquad\qquad\qquad \left. +\Lambda_0(\tau_i^u)e^{f(\vecx_i)}e^{-\Lambda_0(\tau_i^u)e^{f(\vecx_i)}}\right).
\label{app:eq:problematique1}
\end{align}
Second order partial derivatives are given by
\begin{equation}
-\frac{1}{N}\sum_{k:\Delta_k=1}\frac{\partial^2}{\partial f_i \partial f_j}\big(-\Lambda_0(\tau_k)e^{-f(\vecx_k)}\big) = \frac{\delta_{ij}}{N}\Lambda_0(\tau_k)e^{f(\vecx_i)}
\end{equation}
and
\begin{align}
&\sum_{k:\Delta_k=0}\frac{\partial^2}{\partial f_i \partial f_j}\log [S(\tau_i^l)-S(\tau_i^u)] = -\delta_{ij}\left(\frac{\partial}{\partial f_i}\log [S(\tau_i^l)-S(\tau_i^u)]\right)^2\nonumber\\
&\qquad\qquad+ \frac{\delta_{ij}}{[S(\tau_i^l)-S(\tau_i^u)]}\left(-\Lambda_0(\tau_i^l)e^{f(\vecx_i)}e^{-\Lambda_0(\tau_i^l)e^{f(\vecx_i)}} + \left(\Lambda_0(\tau_i^l)e^{f(\vecx_i)}\right)^2e^{-\Lambda_0(\tau_i^l)e^{f(\vecx_i)}}\right.\nonumber\\
&\qquad\qquad\qquad\qquad\left.+\Lambda_0(\tau_i^u)e^{f(\vecx_i)}e^{-\Lambda_0(\tau_i^u)e^{f(\vecx_i)}} - \left(\Lambda_0(\tau_i^u)e^{f(\vecx_i)}\right)^2e^{-\Lambda_0(\tau_i^u)e^{f(\vecx_i)}}\right).
\label{app:eq:problematique2}
\end{align}
Note that (\ref{app:eq:problematique1}) and (\ref{app:eq:problematique2}) can be problematic numerically and the approximations discussed in Section \ref{gphr_num} are required.

\subsection{GP regression with competing risks}
\label{app:mutli}

The negative log likelihood function (with two risk and independent right censoring) is
\begin{align}
\mathcal{L}(\vecf) &= -\frac{1}{N}\sum_{i:\Delta_i\neq1}\log S(t_i|f_i^1) -\frac{1}{N}\sum_{i:\Delta_i\neq2}\log S(t_i|f_i^2) -\frac{1}{N}\sum_{i:\Delta_i=1}\log p(t_i|f_i^1)\nonumber\\
& \quad\quad-\frac{1}{N}\sum_{i:\Delta_i=2}\log p(t_i|f_i^2) + \frac{1}{2N}(\vecf-\ev)\cdot\matK^{-1}(\vecf-\ev) + \log2\pi +\frac{1}{2N}\log|\matK|.
\end{align}
The first order partial derivatives are
\begin{align}
\frac{\partial}{\partial f_i^r} \mathcal{L}(\vecf) &= -\frac{1}{N}\sum_{k:\Delta_k\neq1}\frac{\partial}{\partial f_i^r}\log S(t_k|f_k^1) -\frac{1}{N}\sum_{k:\Delta_k\neq2}\frac{\partial}{\partial f_i^r}\log S(t_k|f_k^2)\nonumber\\
& \quad\quad-\frac{1}{N}\sum_{k:\Delta_k=1}\frac{\partial}{\partial f_i^r}\log p(t_k|f_k^1)-\frac{1}{N}\sum_{k:\Delta_k=2}\frac{\partial}{\partial f_i^r}\log p(t_k|f_k^2)\nonumber\\
&\quad\quad+ \frac{1}{N}[\matK^{-1}(\vecf-\ev)]_i
\end{align}
where 
\begin{equation}
\frac{\partial}{\partial f_i^r}\log p(t_k|f^q_k) = \delta_{ik}\delta_{pq}\beta_q^{-2}(t_k-f_k^q)
\end{equation}
and
\begin{align}
\frac{\partial}{\partial f_i^r}S(t_k|f_k^q) & 
=\delta_{ik}\delta_{pq}\frac{1}{S(t_k|f_k^q)}\int_{t_k}^{\infty}\text{d}s\,\frac{\partial}{\partial f_i^r}\frac{e^{-\frac{1}{2\beta_q^2}(s - f_k^q)^2}}{\sqrt{2\pi}\beta_q}\nonumber\\
& =\delta_{ik}\delta_{pq}\frac{1}{S(t_k|f_k^q)}\frac{1}{\sqrt{2\pi}\beta_q^3}\int_{t_i}^{\infty}\text{d}s\,(s-f_k^q)e^{-\frac{1}{2\beta_q^2}(s - f_k^q)^2}\nonumber\\
&= \delta_{ik}\delta_{pq}\frac{1}{S(t_k|f_k^q)}\frac{1}{\sqrt{2\pi}\beta_q} e^{-\frac{1}{2\beta_q^2}(t_k-f_k^q)^2}.
\end{align}
Second order partial derivatives are (where $\partial^2/\partial f_j^r\partial f_i^q = 0$ for $r\neq q$)
\begin{equation}
\frac{\partial^2}{\partial f_j^r\partial f_i^r}\mathcal{L}(\vecf) =\frac{1}{N}(\matW + \matK^{-1})_{ij}
\end{equation}
with
\begin{align}
\matW_{ij} &= -\sum_{k:\Delta_k\neq1}\frac{\partial^2}{\partial f_i^r\partial f_j^r}\log S(t_k|f_k^1) -\sum_{k:\Delta_k\neq2}\frac{\partial^2}{\partial f_i^r\partial f_j^r}\log S(t_k|f_k^2)\nonumber\\
&-\sum_{k:\Delta_k=1}\frac{\partial^2}{\partial f_i^r\partial f_j^r}\log p(t_k|f_k^1)-\sum_{k:\Delta_k=2}\frac{\partial^2}{\partial f_i^r\partial f_j^r}\log p(t_k|f_k^2).
\end{align}
The matrix $\matW$ is diagonal since
\begin{align}
\frac{\partial^2}{\partial f_i^r\partial f_j^r}\log S(t_k|f_k^q) &= \delta_{ik}\delta_{jk}\left(\frac{1}{S(t_k|f_k^q)}\frac{1}{(2\pi\beta_q^2)^{1/2}} e^{-\frac{1}{2\beta_q^2}(t_k-f_k^q)^2}\right)^2\nonumber\\
& - \frac{(t_k-f_k^q)}{\beta_q^2}\left(\frac{1}{S(t_k|f_k^q)}\frac{1}{(2\pi\beta_q^2)^{1/2}} e^{-\frac{1}{2\beta_q^2}(t_k-f_k^q)^2}\right)
\end{align}
and
\begin{equation}
\frac{\partial^2}{\partial f_i^r\partial f_j^r}\log p(t_k|f_k^q) = \delta_{ik}\delta_{jk}\beta_q^{-2}.
\end{equation}

\bibliographystyle{unsrt}
\bibliography{refs}

\end{document}